\documentclass[11pt,a4paper]{amsart}

\usepackage[latin1]{inputenc}
\usepackage{amsmath}
\usepackage{amsfonts}
\usepackage{amssymb}
\usepackage{amsthm}
\usepackage{hyperref}
\usepackage{anysize}
\marginsize{3.0cm}{3.0cm}{3.0cm}{3.0cm}

\newtheorem{theorem}{Theorem}[section]
\newtheorem{lemma}[theorem]{Lemma}
\newtheorem{proposition}[theorem]{Proposition}
\newtheorem{corollary}[theorem]{Corollary}
\newtheorem{conjecture}[theorem]{Conjecture}

\theoremstyle{definition}
\newtheorem{definition}[theorem]{Definition}
\newtheorem{example}[theorem]{Example}

\newtheorem{question}[theorem]{Question}

\theoremstyle{remark}
\newtheorem{remark}[theorem]{Remark}

\newcommand{\coloneqq}{\mathrel{\mathop:}=}
\newcommand{\form}[1]{{\langle #1 \rangle }}
\newcommand{\pfister}[1]{{\langle \! \langle #1 \rangle \! \rangle}}

\numberwithin{equation}{section}
\begin{document}

\title{Rational Morphisms between Quasilinear Hypersurfaces}
\author{Stephen Scully}
\address{School of Mathematical Sciences, University of Nottingham, University Park, Nottingham, NG7 2RD, United Kingdom}
\email{pmxss4@nottingham.ac.uk}


\begin{abstract} We prove analogues of several well-known results concerning rational morphisms between quadrics for the class of so-called quasilinear $p$-hypersurfaces. These hypersurfaces are nowhere smooth over the base field, so many of the geometric methods which have been successfully applied to the study of projective homogeneous varieties over fields cannot be used. We are therefore forced to take an alternative approach, which is partly facilitated by the appearance of several non-traditional features in the study of these objects from an algebraic perspective. Our main results were previously known for the class of quasilinear quadrics. We provide new proofs here, because the original proofs do not immediately generalise for quasilinear hypersurfaces of higher degree. \end{abstract}

\maketitle

\section{Introduction}

Considerable progress has been made in recent years towards understanding the conditions under which a rational morphism $X \dashrightarrow Y$ can exist between smooth projective varieties $X$ and $Y$ over a field. We can mention, for example, M. Rost's degree formulas and their generalisations (see \cite{Rost}, \cite{Merkurjev}, \cite{LevineMorel}). The results obtained have proven to be particularly successful in their application to the study of projective homogeneous varieties over fields. In the case where $X$ and $Y$ are quadrics, several applications of this approach were already known quite some time ago. Indeed, it was already apparent from the foundational works of A. Pfister in the 1960's that the study of rational morphisms between quadrics has many important applications to the algebraic theory of quadratic forms; for example, to the study of the structure of the Witt ring, or to the construction of fields which exhibit certain arithmetic properties pertaining to quadratic forms. On the other hand, the rich structure theory of the Witt ring unearthed by Pfister permitted the study of rational morphisms between quadrics from a entirely algebraic point of view. In the subsequent decades, the algebraic methods were developed and applied intensively to this area of research. Two of the highlights of this approach are the following results, due to D. Hoffmann and O. Izhboldin respectively:

\begin{theorem}[\cite{Hoffmann1}, Theorem 1] Let $X$ and $Y$ be anisotropic projective quadrics over a field $k$ of characteristic different from 2. If there exists $n \geq 1$ such that $\mathrm{dim}(Y) \leq 2^n - 2 < \mathrm{dim}(X)$, then there are no rational morphisms $X \dashrightarrow Y$. \end{theorem}

\begin{theorem}[\cite{Izhboldin}, Theorem 0.2] Let $X$ and $Y$ be anisotropic projective quadrics over a field $k$ of characteristic different from 2 with $\mathrm{dim}(Y) = 2^n - 1 \leq \mathrm{dim}(X)$ for some $n \geq 0$. If there exists a rational morphism $X \dashrightarrow Y$, then there also exists a rational morphism $Y \dashrightarrow X$. \end{theorem}

Of course, phenomena of this sort are more readily explained nowadays in light of the recent advances in the geometric theory (for example, both the above results can be deduced from Rost's degree formula). Perhaps the most general result available concerning rational morphisms between quadrics is due to N. Karpenko and A. Merkurjev. Before we state it, let us recall a couple of standard definitions. If $X$ is an anisotropic projective quadric over a field $k$ of characteristic different from 2, then the \emph{first Witt index} of $X$, denoted $\mathfrak{i}_1(X)$, is the largest positive integer $i$ such that there exists a degree 1 cycle of dimension $i-1$ on the variety $X_{k(X)}$, where $k(X)$ is the function field of $X$. The \emph{Izhboldin dimension} of $X$, denoted $\mathrm{dim}_{Izh}(X)$, is then defined to be the integer $\mathrm{dim}(X) - \mathfrak{i}_1(X) + 1$. The theorem of Karpenko and Merkurjev now says:

\begin{theorem}[\cite{KarpenkoMerkurjev}, Theorem 4.1] Let $X$ and $Y$ be anisotropic projective quadrics over a field $k$ of characteristic different from 2. If there exists a rational morphism $X \dashrightarrow Y$, then
\begin{enumerate} \item[$\mathrm{(1)}$] $\mathrm{dim}_{Izh}(X) \leq \mathrm{dim}_{Izh}(Y)$.
\item[$\mathrm{(2)}$] $\mathrm{dim}_{Izh}(X) = \mathrm{dim}_{Izh}(Y)$ if and only if there exists a rational morphism $Y \dashrightarrow X$. \end{enumerate} \end{theorem}

If one knows something additional about the first Witt indices of the quadrics involved, then one can start to produce more explicit examples. For instance, Theorems 1.1 and 1.2 can both be recovered from Theorem 1.3 modulo the observation that anisotropic quadrics of dimension $2^n-1$ (for some $n \geq 0$) have first Witt index equal to 1. Although the latter fact was originally proved by D. Hoffmann as a corollary of Theorem 1.1, there are several alternative proofs which are completely independent of Theorem 1.1 (for example, see the paper \cite{Karpenko} of N. Karpenko, where all possible values of the first Witt index are determined). \\

All the above statements include the assumption that the characteristic of the base field is different from 2 (and historically, this is the form in which they were first stated). But it turns out that all these results can be extended to allow for fields of characteristic 2, even when the quadrics involved are not smooth over the base field (see the articles \cite{HoffmannLaghribi1}, \cite{HoffmannLaghribi2} and \cite{Totaro1}). This includes the extreme case of \emph{quasilinear} quadrics, which have no smooth points at all (quasilinear quadrics are those quadrics which are defined by the diagonal part of a symmetric bilinear form over a field of characteristic 2). In fact, several problems which remain open for smooth quadrics have recently been settled for quasilinear quadrics. In particular, we want to mention the results of B. Totaro on the birational geometry of quasilinear quadrics. In the article \cite{Totaro1}, Totaro gives a positive answer to the ``Quadratic Zariski Problem'' for quasilinear quadrics: If $X$ and $Y$ are anisotropic quasilinear quadrics of the same dimension over a field such that there exist rational morphisms from $X$ to $Y$ and from $Y$ to $X$, then $X$ and $Y$ are birational. The proof of this result uses the quasilinear analogue of Theorem 1.3 (due to Totaro) together with the following ``Ruledness Theorem'', proved by Totaro in the same article: If $X$ is an anisotropic quasilinear quadric over a field $k$, then $X$ is birational to $X' \times \mathbb{P}_k^{\mathfrak{i}_1(X) - 1}$ for any subquadric $X' \subset X$ of codimension $\mathfrak{i}_1(X) - 1$ (where the integer $\mathfrak{i}_1(X)$ is defined in the same way as it is for smooth quadrics). For non-quasilinear quadrics, the corresponding problems are still wide open, even in the smooth case.

The goal of the present article is to prove analogues of all the above results for the class of so-called quasilinear $p$-hypersurfaces, which generalises the class of quasilinear quadrics. More precisely, a \emph{quasilinear p-hypersurface} is the projective hypersurface defined by a diagonal form of degree $p$ over a field of characteristic $p>0$. Forms of the latter type are called \emph{quasilinear p-forms}. For evident reasons, the geometry of quasilinear hypersurfaces is only interesting over non-perfect fields. Due in part to the absence of any smooth points on these varieties, one does not have access to the general geometric methods which have proved to be very successful in the study of projective homogeneous varieties over fields. On the other hand, there are several unusual aspects of the theory of quasilinear $p$-forms which suggest that a more algebraic approach is feasible. The main result of the paper is Theorem 5.12. This is an analogue of the main result of \cite{KarpenkoMerkurjev}. The latter result, due to N. Karpenko and A. Merkurjev, immediately implies their Theorem 1.3 as a corollary. Our result should also imply the analogue of Theorem 1.3 for arbitrary quasilinear $p$-hypersurfaces, but there is a technical obstruction which we have currently only managed to overcome for quasilinear quadrics and cubics. This issue is discussed in \S 4. Nevertheless, Theorem 5.12 is already sufficient for several interesting applications which hold for all primes $p$. For example, in \S 6 we prove analogues of Theorems 1.1 and 1.2, as well as some results previously obtained for smooth quadrics by A. Vishik. In the final section, we consider the extension of B. Totaro's results on the birational geometry of quasilinear quadrics to the whole class of quasilinear hypersurfaces. We show that Totaro's ``Ruledness Theorem'' for quasilinear quadrics extends to all primes $p$, but due to the technical issue discussed in \S 4, we are only able to present partial results on the analogue of the ``Quadratic Zariski Problem''. Although all these results were previously known in the case $p=2$, the original proofs do not immediately generalise for larger primes $p$, and we are forced to take a different approach. In particular, we provide new proofs of all the above results for the class of quasilinear quadrics.

The motivation for this paper is the article \cite{Hoffmann2} by D. Hoffmann, where an extensive study of quasilinear $p$-forms was carried out from an algebraic point of view. It was shown there that several aspects of the classical theory of quadratic forms over fields carry over in full generality to this new situation. We recall several notions and results from Hoffmann's paper in \S 2 and \S 3.\\

Throughout this article, $F$ will be a non-perfect field of characteristic $p>0$. $\overline{F}$ will denote a fixed algebraic closure of $F$. By a \emph{scheme} we mean a scheme of finite type over a field. By a \emph{variety}, we mean an integral scheme. A scheme will be called \emph{complete} if it is proper over the base field. If $X$ is a scheme defined over a field $k$, the residue field at a point $x \in X$ will be denoted by $k(x)$. If $X$ is a variety, then we will write $k(X)$ for the function field of $X$. If $L/k$ is any field extension, $X_L$ will denote the scheme $X \times_k \mathrm{Spec}\;L$. Finally, we implicitly assume that all morphisms and rational morphisms of schemes are defined relative to the given base field. 

\section{Quasilinear $p$-forms}

In this section we introduce quasilinear $p$-forms and discuss some of their basic properties and invariants. Everything in this section can be found in the article \cite{Hoffmann2}. Since the proofs of the statements we need are all very short, we include them for the reader's convenience.

\begin{definition} Let $V$ be a finite dimensional $F$-vector space. A \emph{quasilinear p-form} (we will often simply say \emph{form} for simplicity) on $V$ is a map $\varphi \colon V \rightarrow F$ satisfying
\begin{enumerate} \item[$\mathrm{(1)}$] $\varphi(\lambda v) = \lambda^p \varphi(v)$ for all $\lambda \in F$ and all $v \in V$.
\item[$\mathrm{(2)}$] $\varphi(v+w) = \varphi(v) + \varphi(w)$ for all $v,w \in V$. \end{enumerate} \end{definition}

We will say that $\varphi$ is a \emph{quasilinear p-form over F} if $\varphi$ is a quasilinear $p$-form on some finite dimensional $F$-vector space. In this case, we denote the underlying space by $V_\varphi$. The \emph{dimension} of $\varphi$ is the dimension of $V_\varphi$, denoted $\mathrm{dim}(\varphi)$. A \emph{morphism} $\varphi \rightarrow \psi$ of forms over $F$ will be an $F$-vector space morphism $V_\varphi \rightarrow V_\psi$ which carries $\varphi$ to $\psi$. If $\varphi$ and $\psi$ are isomorphic, we will write $\varphi \simeq \psi$. We will say that $\varphi$ \emph{is proportional to $\psi$} if there exists $\lambda \in F^*$ such that $\varphi \simeq \lambda \psi$, where $\lambda \psi$ is the form on $V_\psi$ defined by $v \mapsto \lambda \psi(v)$. If $\varphi$ is a quasilinear $p$-form over $F$, then a \emph{subform} $\psi$ of $\varphi$ is the restriction of $\varphi$ to a subspace of $V_\varphi$. We write $\psi \subset \varphi$. The \emph{direct sum} $\varphi \oplus \psi$ and \emph{tensor product} $\varphi \otimes \psi$ of forms $\varphi, \psi$ over $F$ are defined in the obvious way. If $L/F$ is a field extension and $\varphi$ a quasilinear $p$-form over $F$, we write $\varphi_L$ for the form over $L$ obtained by the extension of scalars. If $\varphi$ is a quasilinear $p$-form over $F$, a vector $v \in V_\varphi$ is called \emph{isotropic} if $\varphi(v) = 0$. We say that the form $\varphi$ is \emph{isotropic} if $V_\varphi$ contains a nonzero isotropic vector; $\varphi$ is called \emph{anisotropic} otherwise. By condition (2) in the definition, the subset of isotropic vectors in $V_\varphi$ is actually a subspace, and $\varphi$ is isotropic if and only if it has nonzero dimension. This additivity condition also implies that $\varphi$ is ``diagonalized" in \emph{every} basis of $V_\varphi$; that is, of the form $a_1x_1^p + ... + a_nx_n^p$. We will use the notation $\form{a_1,...,a_n}$ to denote the quasilinear $p$-form $a_1x_1^p + ... + a_nx_n^p$ on the $F$-vector space $F^n$ in its standard basis. If $\varphi$ is a quasilinear $p$-form over $F$, and $L/F$ is a field extension, we define the \emph{value set}
\begin{equation*} D_L(\varphi) \coloneqq \lbrace \varphi_L(v)\;|\;v \in V_\varphi \otimes_F L \rbrace. \end{equation*}
Since $\varphi$ is additive, $D_L(\varphi)$ is actually an $L^p$-vector subspace of $L$. Clearly it is finite dimensional of dimension $\leq \mathrm{dim}(\varphi)$, and we have $D_L(\varphi) = D_L(\varphi_L)$. \\

The classification of quasilinear $p$-forms over $F$ is given by the following statement:

\begin{proposition}[\cite{Hoffmann2}, Proposition 2.6] Let $\varphi$ be a quasilinear $p$-form over $F$ of dimension $n$, and let $a_1,...,a_m$ be a basis for $D_F(\varphi)$ over $F^p$. Then $m \leq n$ and there is an isomorphism of forms over $F$:
\begin{equation*} \varphi \simeq \form{a_1,...,a_m} \oplus \form{\underbrace{0,...,0}_{n-m}}. \end{equation*}
\begin{proof} Let $W \subseteq V_\varphi$ be the subspace of isotropic vectors. For each $1 \leq i \leq m$, let $v_i \in V_\varphi$ be such that $\varphi(v_i) = a_i$, and let $U$ be the subspace of $V_\varphi$ spanned by the $v_i$. Since the $a_i$ are linearly independent over $F^p$, the $v_i$ are linearly independent over $F$. Therefore, in order to prove the statement, it suffices to show that $V_\varphi = U \oplus W$. Clearly $U \cap W = \lbrace 0 \rbrace$. On the other hand, given any $v \in V_\varphi$, we can find $\lambda_i \in F$ such that $\varphi(v) = \sum_{i=1}^m \lambda_i^p a_i$. Then $\varphi(v) = \varphi(\sum_{i=1}^m \lambda_iv_i)$, so that $v - \sum_{i=1}^m \lambda_i v_i$ is an isotropic vector. Therefore $v \in U + W$, and the statement is proved. \end{proof} \end{proposition}

The isomorphism class of a quasilinear $p$-form $\varphi$ over $F$ is therefore determined by two invariants: the dimension of the subspace of isotropic vectors, and the $F^p$-vector space $D_F(\varphi)$. In particular, the theory is essentially vacuous over perfect fields (i.e. when $F = F^p$), which is why we assume $F$ to be non-perfect. In the case of anisotropic forms, we get:

\begin{corollary} Let $\varphi$ and $\psi$ be anisotropic quasilinear $p$-forms over $F$. Then $\psi$ is isomorphic to a subform of $\varphi$ if and only if $D_F(\psi) \subset D_F(\varphi)$. \end{corollary}

If $\varphi$ is an arbitrary quasilinear $p$-form over $F$, it follows from Proposition 2.2. that there is a unique (up to isomorphism) anisotropic subform $\varphi_{an} \subset \varphi$ such that $\varphi = \varphi_{an} \oplus \form{0,...,0}$. The form $\varphi_{an}$ is called the \emph{anisotropic part} of $\varphi$. The integer $\mathfrak{i}_0(\varphi) \coloneqq \mathrm{dim}(\varphi) - \mathrm{dim}(\varphi_{an})$ is called the \emph{defect index} of $\varphi$. By the proof of Proposition 2.2, $\mathfrak{i}_0(\varphi)$ is nothing else but the dimension of the subspace of isotropic vectors in $V_\varphi$. We write $\varphi \sim \psi$ whenever $\varphi_{an} \simeq \psi_{an}$ for $p$-forms $\varphi$, $\psi$ over $F$. This is obviously an equivalence relation on the set of isomorphism classes of quasilinear $p$-forms over $F$. The following consequence of Proposition 2.2 will be used several times in the sequel:

\begin{lemma} Let $\varphi = \form{a_1,...,a_n}$ be a quasilinear $p$-form over $F$, and let $L/F$ be any field extension. Then there is a subset $\lbrace a_{j_1},...,a_{j_m} \rbrace \subset \lbrace a_1,...,a_n \rbrace$ such that $(\varphi_L)_{an} \simeq \form{a_{j_1},...,a_{j_m}}$.
\begin{proof} Since $D_L(\varphi) = L^p \cdot a_1 + ... + L^p \cdot a_n$, we can find a subset $\lbrace a_{j_1},...,a_{j_m} \rbrace \subset \lbrace a_1,...,a_n \rbrace$ which constitutes a basis of $D_L(\varphi)$ over $L^p$. The statement then follows from Proposition 2.2. \end{proof} \end{lemma}

We now introduce a special class of forms. In the theory of quadratic forms over fields of characteristic $\neq 2$, an important role is played by the class of so-called \emph{Pfister forms} (namely, tensor products of binary forms which represent 1). The set of isomorphism classes of $n$-fold Pfister forms over such a field $k$ is in bijection with the set of pure symbols in the torsion Milnor $K$-group $K_n^M(k)/2K_n^M(k)$. The projective quadric defined by an $n$-fold Pfister form $\pfister{a_1,...,a_n}$ is a \emph{splitting variety} in the sense that for any extension $L/k$, the quadric has an $L$-rational point if and only if the symbol $\lbrace a_1,...,a_n \rbrace$ is divisible by 2 in $K_n^M(L)$. Let us now recall the equal characteristic analogue of the norm residue isomorphism theorem (formerly the Bloch-Kato conjecture), due to S. Bloch, K. Kato and O. Gabber:

\begin{theorem} [\cite{Kato}, \cite{BlochKato}] For any field $F$ of characteristic $p>0$, there are isomorphisms of abelian groups
\begin{eqnarray*} K_n^M(F)/pK_n^M(F) &\xrightarrow{\sim}& \mathrm{ker}(\wp \colon \Omega_F^n \rightarrow \Omega_F^n/d\Omega_F^{n-1})\\
\lbrace a_1,...,a_n \rbrace & \mapsto & \frac{da_1}{a_1} \wedge ... \wedge \frac{da_n}{a_n}. \end{eqnarray*}
\end{theorem}

Here, $(\Omega_F^\bullet, d)$ is the de Rham complex of absolute differential forms over $F$, and $\wp$ is the inverse Cartier operator (see \cite{Kato} for details). It follows from this result that for any $a_1,...,a_n \in F^*$, the pure symbol $\lbrace a_1,...,a_n \rbrace$ is divisible by $p$ in $K_n^M(F)$ if and only if $[F^p(a_1,...,a_n):F^p] <p^n$. Consider the quasilinear $p$-form
\begin{equation*} \pfister{a_1,...,a_n} \coloneqq \sum_{0 \leq j_1,...,j_n \leq p-1}(\prod_{i=1}^n a_i^{j_i})x_{j_1,...,j_n}^p \end{equation*}
of dimension $p^n$ on the $F$-vector space $F^{p^n}$ in its standard basis. Clearly $\pfister{a_1,...,a_n} \simeq \pfister{a_1} \otimes ... \otimes \pfister{a_n}$. By the definition, we have $D_L(\pfister{a_1,...,a_n}) = L^p(a_1,...,a_n)$ for every field extension $L/F$. It therefore follows from Proposition 2.2 that the form $\pfister{a_1,...,a_n}_L$ is isotropic if and only if $[L^p(a_1,...,a_n):L^p] < p^n$, which in turn holds if and only if $\lbrace a_1,...,a_n \rbrace = 0$ in $K_n^M(L)/pK_n^M(L)$. The (projective) hypersurfaces defined by the forms $\pfister{a_1,...,a_n}$, which are called \emph{quasi-Pfister forms}, may therefore be regarded (in the above sense) as splitting varieties in the equal characteristic setting. Moreover, it turns out that the quasi-Pfister forms exhibit all the same properties as the Pfister quadratic forms in the mixed characteristic: they are precisely those forms which are ``multiplicative''; up to proportionality, the anisotropic quasi-Pfister forms are characterised by the degree to which they split over the generic point of the associated hypersurface; they have the ``roundness'' property, etc. We refer to the article \cite{Hoffmann2} for further details.\\

It will be convenient to record here the following observation concerning the splitting of quasi-Pfister forms over extensions of the base field:

\begin{lemma} Let $\pi = \pfister{a_1,...,a_n}$ be an anisotropic quasi-Pfister form over $F$. Let $L/F$ be a field extension such that $\pi_L$ is isotropic. Then there is a proper subset $\lbrace a_{j_1},...,a_{j_m} \rbrace \subset \lbrace a_1,...,a_n \rbrace$ such that $(\pi_L)_{an} \simeq \pfister{a_{j_1},...,a_{j_m}}$.
\begin{proof} Let $m$ be such that $[L^p(a_1,...,a_n):L^p] = p^m$. Since $\pi_L$ is isotropic, we have $m<n$. Now, we may choose a subset $\lbrace a_{j_1},...,a_{j_m} \rbrace$ of  $\lbrace a_1,...,a_n \rbrace$ such that $L^p(a_1,...,a_n) = L^p(a_{j_1},...,a_{j_m})$. Then the quasi-Pfister form $\pfister{a_{j_1},...,a_{j_m}}$ is anisotropic over $L$, and we have
\begin{equation*} D_L((\pi_L)_{an}) = D_L(\pi) = L^p(a_1,...,a_n) = L^p(a_{j_1},...,a_{j_m}) = D_L(\pfister{a_{j_1},...,a_{j_m}}). \end{equation*}
The statement therefore follows from Corollary 2.3. \end{proof} \end{lemma}

To the arbitrary form $\varphi$, we associate a certain anisotropic quasi-Pfister form. Let us first recall the following invariant:

\begin{definition}[\cite{HoffmannLaghribi1}, \cite{Hoffmann2}] Let $\varphi$ be a quasilinear $p$-form over a field $F$, and let $L/F$ be a field extension. The \emph{norm field} of $\varphi$ (over $L$) is the field
\begin{equation*} N_L(\varphi) \coloneqq L^p(\frac{a}{b}\;|\;a,b \in D_L(\varphi) \cap L^*). \end{equation*} \end{definition}

Note that if two forms $\varphi$ and $\psi$ are proportional over $F$, then $N_L(\varphi) = N_L(\psi)$ for all extensions $L/F$. The norm field invariant was first introduced in the context of quasilinear quadratic forms by D. Hoffmann and A. Laghribi in \cite{HoffmannLaghribi1}. It turns out that this is a birational invariant of quasilinear $p$-hypersurfaces (see Proposition 4.10), and we will exploit this for the proofs of our main results. A more direct description of the norm field is given by the following lemma:

\begin{lemma} Let $\varphi$ denote the quasilinear $p$-form $\form{a_0,...,a_n}$ over $F$. Then
\begin{equation*} N_F(\varphi) = F^p(\frac{a_1}{a_0},...,\frac{a_n}{a_0}) \end{equation*}
\begin{proof} Since the norm field does not change when we multiply $\varphi$ by a scalar, we may assume that $a_0 = 1$. In this case, it is clear that $F^p(a_1,...,a_n) \subseteq N_F(\varphi)$. But the reverse inclusion also holds, since $D_F(\varphi) = F^p + F^p \cdot a_1 + ... + F^p \cdot a_n \subseteq F^p(a_1,...,a_n)$. \end{proof} \end{lemma}

It follows that $N_L(\varphi)$ is finite dimensional over its subfield $L^p$ for any form $\varphi$ over $F$ and any extension $L/F$. Moreover, if $a_1,...,a_m \in F$ are such that $N_F(\varphi) = F^p(a_1,...,a_m)$, then we have $N_L(\varphi) = N_L(\varphi_L) = L^p(a_1,...,a_m)$ for all extensions $L/F$. Note that the dimension of $N_L(\varphi)$ as an $L^p$-vector space is always a power of $p$. One may therefore define:

\begin{definition} The integer $\mathrm{lndeg}_L(\varphi) \coloneqq \mathrm{log}_p([N_L(\varphi):L^p])$ is called the (\emph{logarithmic}) \emph{norm degree} of $\varphi$ (over $L$). \end{definition}

\begin{remark} In \cite{Hoffmann2}, the \emph{norm degree} of a form $\varphi$  over $F$ is defined to be the integer $\mathrm{ndeg}_F(\varphi) \coloneqq p^{\mathrm{lndeg}_F(\varphi)}$. For our purposes, it will be more convenient to take the base $p$-logarithm; for example, see Theorem 4.2 for a result of S. Schr\"{o}er which shows that the integer $\mathrm{lndeg}_F(\varphi)$ determines the size of the singular (non-regular) locus of the projective hypersurface $\lbrace \varphi = 0 \rbrace$. \end{remark}

Finally, we define:

\begin{definition}[\cite{HoffmannLaghribi1}, \cite{Hoffmann2}] Let $\varphi$ be a $p$-form over $F$, and let $L/F$ be a field extension. If $\mathrm{lndeg}_L(\varphi) = m$, and $a_1,...,a_m \in L^*$ are such that $N_L(\varphi) = L^p(a_1,...,a_m)$, then the anisotropic quasi-Pfister form $\widehat{\nu}_L(\varphi) \coloneqq \pfister{a_1,...,a_m}$ is called the \emph{norm form} of $\varphi$ (over $L$). \end{definition}

By Corollary 2.3, this does not depend on the choice of generators $a_i$ up to isomorphism. Note that we have $\widehat{\nu}_L(\varphi) = \widehat{\nu}_L((\varphi_L)_{an})$, and by the proof of Lemma 2.6, both forms are equal to the anisotropic part of $\widehat{\nu}_F(\varphi)$ over $L$. By definition, the dimension of $\widehat{\nu}_L(\varphi)$ is nothing else but $p^{\mathrm{lndeg}_L(\varphi)}$. If $\varphi_L$ is anisotropic, then by Corollary 2.3 and Lemma 2.8, $\varphi_L$ is proportional to a subform of $\widehat{\nu}_L(\varphi)$.

\section{Quasilinear $p$-forms over extensions of the base field}

In this section we record some general observations about the behaviour of quasilinear $p$-forms and their invariants over extensions of the base field. Again, most of the statements here can be found in the article \cite{Hoffmann2}.\\

The invariant $\mathrm{lndeg}_F(\varphi)$ defined above gives a useful necessary (but not sufficient) condition for an anisotropic quasilinear $p$-form to become isotropic over an extension of the base field.

\begin{proposition}[\cite{Hoffmann2}, Proposition 5.2] Let $\varphi$ be an anisotropic quasilinear $p$-form over $F$, and let $L/F$ be a field extension over which $\varphi$ becomes isotropic. Then $\mathrm{lndeg}_L(\varphi) < \mathrm{lndeg}_F(\varphi)$.
\begin{proof} Recall that for any extension $E/F$, the integer $p^{\mathrm{lndeg}_E(\varphi)}$ is equal to the dimension of the anisotropic part of the norm form $\widehat{\nu}_F(\varphi)$ over $E$. Proving the statement therefore amounts to checking that $\widehat{\nu}_F(\varphi)$ becomes isotropic over $L$. But this is evident, because $\varphi$ is proportional to a subform of $\widehat{\nu}_F(\varphi)$. \end{proof} \end{proposition}

\begin{corollary} Let $\varphi$ be an anisotropic quasilinear $p$-form over $F$, and let $a \in F \setminus F^p$. If $\varphi_{F(\sqrt[p]{a})}$ is isotropic, then $a \in N_F(\varphi)$.
\begin{proof} Let $m = \mathrm{lndeg}_F(\varphi)$, and let $a_1,...,a_m \in F^*$ be such that $N_F(\varphi) = F^p(a_1,...,a_m)$. By Proposition 3.1, the field $N_{F(\sqrt[p]{a})}(\varphi) = F(\sqrt[p]{a})^p(a_1,...,a_m) = F^p(a_1,...,a_m,a)$ has dimension $\leq p^{m-1}$ over $F^p(a) = F(\sqrt[p]{a})^p$. It must therefore have dimension $p^m$ over $F^p$, and so $a \in F^p(a_1,...,a_m) = N_F(\varphi)$. \end{proof} \end{corollary}

This allows us to give another characterisation of the norm field for anisotropic forms:

\begin{corollary} Let $\varphi$ be an anisotropic quasilinear $p$-form over $F$. Let $S$ be the set of all $a \in F$ such that $\varphi_{F(\sqrt[p]{a})}$ is isotropic. Then $N_L(\varphi) = L^p(S)$ for all field extensions $L/F$. \begin{proof} By the remarks following Lemma 2.8, it suffices to prove this for $L=F$. Since $N_F(\varphi)$ is invariant under multiplying $\varphi$ by a scalar, we may assume that $\varphi = \form{1,a_1,...,a_n}$ for some $a_i \in F^*$. Then $N_F(\varphi) = F^p(a_1,...,a_n)$ by Lemma 2.8. Since $\varphi_{F(\sqrt[p]{a_i})}$ is clearly isotropic for all $1 \leq i \leq n$, we have $N_F(\varphi) \subset F^p(S)$. The reverse inclusion follows from Corollary 3.2. \end{proof} \end{corollary}

Let $L/F$ be a finitely generated field extension. Recall that $L/F$ is called \emph{separably generated} if $L$ can be realised as a purely transcendental extension of $F$ followed by a separable algebraic extension. If $L/F$ is separably generated, then there are several rather direct ways to see that any anisotropic quasilinear $p$-form over $F$ remains anisotropic over $L$. More generally, we have:

\begin{proposition} Let $L/F$ be any field extension. Then there exists an anisotropic quasilinear $p$-form $\varphi$ over $F$ such that $\varphi_L$ is isotropic if and only if $L/F$ is not separably generated.
\begin{proof} The extension $L/F$ is not separably generated if and only if $L^p/F^p$ is not separably generated. By Proposition 4.1 in Chapter VIII of \cite{Lang} (we refer to the proof rather than the statement itself), the latter is true if and only if $L^p$ and $F$ are not linearly disjoint over $F^p$, which is a precise translation of the condition in our statement. \end{proof}  \end{proposition}

In particular, we have:

\begin{corollary}[\cite{Hoffmann2}, Proposition 5.3] Let $\varphi$ be an anisotropic quasilinear $p$-form over $F$, and let $L/F$ be a separably generated field extension. Then $\varphi_L$ is anisotropic. Moreover, $\mathrm{lndeg}_L(\varphi) = \mathrm{lndeg}_F(\varphi)$.\begin{proof} The first assertion is implicit in Proposition 3.4. The second statement follows by applying the first statement to the norm form $\widehat{\nu}_F(\varphi)$. \end{proof} \end{corollary}

\begin{remark} For finite separable extensions, the first part of this statement may be viewed as a replacement for Springer's Theorem from the theory of quadratic forms, which says that an anisotropic quadratic form remains anisotropic over any odd degree extension of the base field. \end{remark}

Thus in order to study the isotropy behaviour of quasilinear $p$-forms over extensions of the base field, we are essentially reduced to considering finite purely inseparable extensions. We conclude this section by collecting a couple of simple facts with this in mind.

\begin{lemma} Let $\varphi$ be a quasilinear $p$-form over $F$, and let $L/F$ be a finite extension of degree $n$. If $\varphi_L$ is isotropic, then $\varphi$ contains a subform of dimension $\leq n$ which becomes isotropic over $L$.
\begin{proof} Let $\mu_1,...,\mu_n$ be a basis for $L$ over $F$, and suppose that $w \in V_\varphi \otimes_F L$ is a nonzero isotropic vector for $\varphi_L$. Then we can write
\begin{equation*} w = v_1 \otimes \mu_1 + ... + v_n \otimes \mu_n \end{equation*}
for some $v_i \in V_\varphi$. The $v_i$ span a nonzero subspace $W \subset V_\varphi$ of dimension $\leq n$, and the restriction $\varphi|_W$ becomes isotropic over $L$. \end{proof} \end{lemma}

\begin{lemma} Let $\varphi$ be an anisotropic quasilinear $p$-form over $F$, and let $L/F$ be a degree $p$ extension. Then 
\begin{enumerate} \item[$\mathrm{(1)}$] $\mathrm{dim}(\varphi_L)_{an} \geq \frac{1}{p} \mathrm{dim}(\varphi)$.
\item[$\mathrm{(2)}$] $\mathrm{lndeg}_L(\varphi) \geq \mathrm{lndeg}_F(\varphi) - 1$, with equality if $\varphi_L$ is isotropic. \end{enumerate}
\begin{proof} For any extension $E/F$, the dimension of ($\varphi_{E})_{an}$ is equal to the dimension of the $E^p$-vector space $D_E(\varphi)$ by Proposition 2.2. In order to prove (1), we therefore have to show that $\mathrm{dim}_{L^p}(D_L(\varphi)) \geq \frac{1}{p}\mathrm{dim}_{F^p}(D_F(\varphi))$. But this is obvious, because
\begin{equation*} p \cdot \mathrm{dim}_{L^p}(D_L(\varphi)) = \mathrm{dim}_{F^p}(D_L(\varphi)) \geq \mathrm{dim}_{F^p}(D_F(\varphi)). \end{equation*}
The inequality in statement (2) now follows by applying (1) to the norm form $\widehat{\nu}_F(\varphi)$. In the case where $\varphi_L$ is isotropic, equality holds by Proposition 3.1. \end{proof} \end{lemma}

\section{Quasilinear $p$-hypersurfaces}

In this section we present some basic observations concerning quasilinear $p$-hypersurfaces and rational morphisms between them. We start with the following lemma:

\begin{lemma} Let $F$ be a field of characteristic $p>0$. Let $f = a_0x_0^p + a_1x_1^p + ... + a_nx_n^p \in F[x_0,...,x_n]$ be a polynomial, and assume that $a_0 \neq 0$. Then $f$ is reducible in $F[x_0,...,x_n]$ if and only if $\frac{a_i}{a_0} \in F^p$ for all $i \in [1,n]$.
\begin{proof} We may assume that $a_0 = 1$. If $a_i \in F^p$ for all $i \in [1,n]$, then $f = (\sum_{i=0}^n \sqrt[p]{a_i} x_i)^p$. Conversely, if $f$ is reducible in $F[x_0,...,x_n]$, then it is certainly reducible in $F(x_1,...,x_n)[x_0]$. It follows that $a_1x_1^p + ... + a_nx_n^p$ is a $p^{th}$-power in $F(x_1,...,x_n)$, and this easily implies that $a_i \in F^p$ for all $i \in [1,n]$. \end{proof} \end{lemma}

Now, let $\varphi$ be a quasilinear $p$-form over $F$ of dimension $d+2$. We consider the projective hypersurface $X_\varphi \coloneqq \mathrm{Proj}(S^\bullet(V_\varphi^*)/(\varphi)) \subset \mathbb{P}(V_\varphi)$ of dimension $d$ over $F$. A scheme of this type will be called a \emph{quasilinear p-hypersurface}. By Lemmas 2.8 and 4.1, the scheme $X_\varphi$ is integral if and only if $\mathrm{lndeg}_F(\varphi) >0$. In particular, $X_\varphi$ is integral whenever $\varphi$ is anisotropic. If $\mathrm{lndeg}_F(\varphi)>0$, then we let $F(\varphi)$ denote the field of rational functions on $X_\varphi$. Note that $F(\varphi)$ may be realised as a purely transcendental extension of $F$ followed by a purely inseparable extension of degree $p$. If $L/F$ is any extension of fields, then $X_{\varphi_L}$ is canonically isomorphic to $(X_\varphi)_L$, and by construction, $X_\varphi$ has an $L$-rational point if and only if $\varphi_L$ is isotropic. In particular, the anisotropic form $\varphi$ becomes isotropic over the field $F(\varphi)$. We will say that $X_\varphi$ is isotropic (resp. anisotropic) if $\varphi$ is isotropic (resp. anisotropic), and we define $\mathfrak{i}_0(X_\varphi) \coloneqq \mathfrak{i}_0(\varphi)$. By Corollary 3.5, $X_\varphi$ is isotropic if and only if it has a zero cycle of degree 1. Moreover, if $X_\varphi$ is isotropic, then Proposition 2.2 shows that $X_\varphi$ is a cone over $X_{\varphi_{an}}$ with vertex given by the linear subspace of dimension $\mathfrak{i}_0(X_\varphi) - 1$ corresponding to the subspace of isotropic vectors in $V_\varphi$. It follows that for any $i \geq 0$, we have $\mathfrak{i}_0(X_\varphi)> i$ if and only if $X_\varphi$ has a dimension $i$ cycle of degree prime to $p$ (where by \emph{degree} we mean the degree as a cycle on the ambient projective space). \\

A quasilinear $p$-hypersurface $X$ over $F$ is a twisted form of the $p^{th}$ infinitesimal neighbourhood of a hyperplane in some projective space $\mathbb{P}^n_{\overline{F}}$. In particular, the smooth locus of $X$ is empty. Still, since the base field $F$ is assumed to be non-perfect, it is interesting to ask when $X$ is a regular scheme. The following result, due to S. Schr\"{o}er, shows that this is rarely the case:

\begin{theorem}[\cite{Schroer}, Theorem 3.3] Let $\varphi$ be a quasilinear $p$-form over $F$ of dimension $\geq 2$. Then $X_\varphi$ is a regular scheme if and only if $\mathrm{lndeg}_F(\varphi) = \mathrm{dim}(\varphi) - 1$. If $X_\varphi$ is not regular, then the non-regular locus has codimension $\mathrm{lndeg}_F(\varphi)$ in $X_\varphi$. \end{theorem}

The remainder of this section consists of some general observations concerning rational morphisms between quasilinear $p$-hypersurfaces. Since an isotropic quasilinear $p$-hypersurface is a cone over its anisotropic part (and hence stably birational to its anisotropic part), we may restrict our attention to the anisotropic case. Note that if $\varphi$ and $\psi$ are quasilinear $p$-forms over $F$ of dimension $\geq 2$ with $\mathrm{lndeg}_F(\psi) > 0$, then the existence of a rational map $X_\psi \dashrightarrow X_\varphi$ is equivalent to the isotropy of the form $\varphi_{F(\psi)}$. Indeed, given a rational map $X_\psi \dashrightarrow X_\varphi$, the closure of its graph in $X_\psi \times X_\varphi$ pulls back to a rational point on the generic fibre $(X_\varphi)_{F(\psi)}$ of the canonical projection $X_\psi \times X_\varphi \rightarrow X_\psi$. Conversely, any rational point of $(X_\varphi)_{F(\psi)}$ can be viewed as the generic point of a closed subvariety of $X_\psi \times X_\varphi$ birational to $X_\psi$ over $F$; using the other projection $X_\psi \times X_\varphi \rightarrow X_\varphi$, we get a rational morphism $X_\psi \dashrightarrow X_\varphi$. We will switch between the algebraic and geometric terminology where we feel it is appropriate.\\

It will be useful to observe the following simple fact:

\begin{lemma} Let $\varphi$ and $\psi$ be anisotropic quasilinear $p$-forms over $F$. Then
\begin{enumerate} \item[$\mathrm{(1)}$] $\mathrm{dim}(\varphi_{F(\psi)})_{an} \geq \frac{1}{p} \mathrm{dim}(\varphi)$.
\item[$\mathrm{(2)}$] $\mathrm{lndeg}_{F(\psi)}(\varphi) \geq \mathrm{lndeg}_F(\varphi) - 1$, with equality if $\varphi_{F(\psi)}$ is isotropic. \end{enumerate}
\begin{proof} Since $F(\psi)$ can be realised as a purely transcendental extension of $F$ followed by a degree $p$ extension, this follows from Corollary 3.5 and Lemma 3.8. \end{proof} \end{lemma}

Now, a basic question here is the following:

\begin{question} Let $X$, $Y$ and $Z$ be anisotropic quasilinear $p$-hypersurfaces over $F$. Suppose that there exist rational morphisms $X \dashrightarrow Y$ and $Y \dashrightarrow Z$. Does there exist a rational morphism $X \dashrightarrow Z$? \end{question}

One approach to this problem is suggested by the following classical result:

\begin{proposition} Let $X,Y$ and $Z$ be varieties over a field $k$ with $Z$ complete, and suppose that there are rational morphisms $X \dashrightarrow Y$ and $Y \dashrightarrow Z$. If the image of the generic point of $X$ under the map $X \dashrightarrow Y$ is a regular point of $Y$, then there exists a rational morphism $X \dashrightarrow Z$.
\begin{proof} This is a standard application of the valuative criterion of properness. More explicitly, let $y \in Y$ be the image of the generic point of $X$ under the map $X \dashrightarrow Y$. Since $y$ is regular, there is a valuation ring $R$ of the function field $k(Y)$ with residue field $k(y)$. By the valuative criterion of properness, the map $\mathrm{Spec}\;k(Y) \rightarrow Z$ extends to a morphism $\mathrm{Spec}\;R \rightarrow Z$. Passing to the residue field we, get a morphism $\mathrm{Spec}\;k(y) \rightarrow Z$. Finally, composing this with the natural map $\mathrm{Spec}\;k(X) \rightarrow \mathrm{Spec}\;k(y)$ gives a morphism $\mathrm{Spec}\;k(X) \rightarrow Z$, as we wanted. \end{proof} \end{proposition}

\begin{corollary} Let $Y$ be a complete variety over $F$, and let $\varphi$ and $\psi$ be anisotropic quasilinear $p$-forms of dimension $\geq 2$ over $F$ such that $\psi$ is proportional to a subform of $\varphi$. If there exists a rational morphism $X_\varphi \dashrightarrow Y$, then there exists a rational morphism $X_\psi \dashrightarrow Y$.
\begin{proof} It suffices to treat the case where $\psi$ is a codimension 1 subform of $\varphi$. In this case, $X_\varphi$ is regular at the generic point of $X_\psi$, because the latter is an effective Cartier divisor in $X_\varphi$. The statement therefore follows from Proposition 4.5. \end{proof} \end{corollary} 

Now observe that the statement of Question 4.4 depends not on the variety $X$, but only on its generic point. Moreover, the function field $F(X)$ can be realised as a purely transcendental extension of $F$ followed by a purely inseparable extension of degree $p$. Replacing the base field $F$ with a suitable purely transcendental extension of it, we therefore reduce to the case where $X$ is just a point (of degree $p$). Using Lemma 3.7 and Corollary 4.6, we can further reduce to the case where $\mathrm{dim}(Y) \leq p-2$. Taking Proposition 4.5 into consideration, we see that Question 4.4 can be settled affirmatively with a positive answer to the following question:

\begin{question} Let $Y$ be an anisotropic quasilinear $p$-hypersurface of dimension $\leq p-2$ over $F$, and let $y \in Y$ be a closed point of degree $p$. Is it true that there exists a regular closed point $z \in Y$ such that $F(z) \cong F(y)$ over $F$? \end{question}

Of course, this is trivially true in the case where $p=2$. Therefore Question 4.4 has a positive answer for quasilinear quadrics, as was well-known previously. For larger primes $p$, it was essentially asked in \cite{Hoffmann2} whether Question 4.7 has a positive answer with the much stronger condition that the point $z$ be the intersection of the hypersurface $Y$ with a line in the ambient projective space (see Question B in \S 5 of that paper). In general, this is not the case, even for $p=3$. To give an explicit example, let $a,b \in F^*$ be such that $[F^3(a,b):F^3] = 3^2 = 9$. Then the anisotropic cubic form $\varphi = \form{1,a+ba^2,b}$ becomes isotropic over the field $F(\sqrt[3]{a})$, but one easily checks that $\varphi$ has no 2-dimensional subforms which become isotropic over the same extension. Nevertheless, the weaker assertion we make here holds in the case where $p=3$:

\begin{proposition} Question 4.7 (and hence Question 4.4) has a positive answer for $p=2$ and $p=3$.
\begin{proof} If $Y$ is just a point, then there is nothing to prove. We may therefore assume that $p=3$ and $\mathrm{dim}(Y) = 1$. Let $\varphi$ be 3-dimensional anisotropic form over $F$ which defines the hypersurface $Y$. By Lemma 2.8, we either have $\mathrm{lndeg}_F(\varphi) = 1$ or $\mathrm{lndeg}_F(\varphi) = 2$. In the latter case, $Y$ is a regular at all of its points by Theorem 4.2, and again there is nothing to prove. We are therefore left with the case where $\mathrm{lndeg}_F(\varphi) = 1$. Since $\varphi$ becomes isotropic over the residue field $F(y)$, it follows from Proposition 3.1 that we have $\mathrm{lndeg}_{F(y)}(\varphi) = 0$. In other words, the subspace of isotropic vectors for $\varphi_{F(y)}$ is 2-dimensional. For dimension reasons, it follows that every 2-dimensional subform of $\varphi$ over $F$ becomes isotropic over $F(y)$. So if $z$ is a closed point of $Y$ defined by any 2-dimensional subform of $\varphi$, we therefore have $F(z) \cong F(y)$ over $F$. Since all such points are regular, the statement is proved. \end{proof} \end{proposition}

Actually, the proof of Proposition 4.8 shows that for \emph{any} prime $p$, Question 4.7 has a positive answer whenever $\mathrm{dim}(Y) \leq 1$. In particular, for the prime 5 we only need to treat the case where $\mathrm{dim}(Y) = 2$ or $\mathrm{dim}(Y) = 3$. Suppose that $\mathrm{dim}(Y) = 2$, and let $\varphi$ be a 4-dimensional anisotropic form defining $Y$. By Lemma 2.8, we have $1 \leq \mathrm{lndeg}_F(\varphi) \leq 3$. Again, the proof of Proposition 4.8 shows that Question 4.7 has a positive answer if $\mathrm{lndeg}_F(\varphi) = 1$ or $\mathrm{lndeg}_F(\varphi) = 3$. Therefore the only interesting case is where $\mathrm{lndeg}_F(\varphi) = 2$. The following example suggests that things are already rather more complicated in this situation:

\begin{example} Let $p=5$, and let $a,b \in F^*$ be such that $[F^5(a,b):F^5] = 5^2 = 25$. Choose a polynomial $g \in F^5[s,t]$ of degree $\leq 4$ in both variables $s,t$ so that the form
\begin{equation*} \varphi = \form{1,a,b,g(a,b)} \end{equation*}
is anisotropic over $F$. By the definition of $\varphi$, we have $\mathrm{lndeg}_F(\varphi) = 2$. Moreover, condition on the coefficients $a,b$ implies that there are derivations $D_a, D_b \colon F \rightarrow F$ satisfying
\begin{enumerate} \item[$\mathrm{(1)}$] $D_a(a) = 1$, $D_a(b) = 0$, and
\item[$\mathrm{(2)}$] $D_b(a) = 0$, $D_b(b) = 1$. \end{enumerate}
For indeterminates $T_0,T_1,T_2,T_3$, we can extend these to derivations $D_a, D_b \colon F[T_0,...,T_3] \rightarrow F[T_0,...,T_3]$ by sending the variables to zero. Now, let us identify our form $\varphi$ with the polynomial $T_0^5 + aT_1^5 + bT_2^5 + g(a,b)T_3^5 \in F[T_0,...,T_3]$. Then the derivatives
\begin{equation*} D_a(\varphi) = T_1^5 + \frac{\partial g}{\partial s}(a,b)T_3^5, \end{equation*}
and
\begin{equation*} D_b(\varphi) = T_2^5 + \frac{\partial g}{\partial t}(a,b)T_3^5 \end{equation*}
necessarily vanish at the non-regular points of the hypersurface $X_\varphi$. A direct calculation then shows that the non-regular locus of $X_\varphi$ consists of just one point, with residue field isomorphic to
\begin{equation*} L = F(\sqrt[5]{\frac{\partial g}{\partial s}(a,b)}, \sqrt[5]{\frac{\partial g}{\partial t}(a,b)}, \sqrt[5]{g(a,b) - a\frac{\partial g}{\partial s}(a,b) - b\frac{\partial g}{\partial t}(a,b)}). \end{equation*}
As far as Question 4.7 is concerned, we are only interested in points of minimal degree (in this case, degree 5), so the first task here is to determine conditions on the polynomial $g$ under which $[L:F] = 5$. One can hope that the restrictions imposed on $g$ are sufficiently strong to force the existence of more than one $L$-rational point on the hypersurface $X_\varphi$; by the proof of Proposition 4.8, this would be enough to provide a positive answer to Question 4.7 for the form $\varphi$. \end{example}

We conclude this section by pointing out the following important consequence of Corollary 4.6. It was first proved by D. Hoffmann and A. Laghribi for quasilinear quadratic forms in \cite{HoffmannLaghribi1}, and was later extended to arbitrary quasilinear $p$-forms by D. Hoffmann in \cite{Hoffmann2}:

\begin{proposition}[\cite{Hoffmann2}, Lemma 7.12] Let $\varphi$ and $\psi$ be anisotropic quasilinear $p$-forms over $F$ of dimension $\geq 2$. If there exists a rational morphism $X_\psi \dashrightarrow X_\varphi$, then $N_F(\psi) \subset N_F(\varphi)$.
\begin{proof} Since any two forms which are proportional have the same norm field, we may assume that $\psi = \form{1,a_1,...,a_n}$ for some $a_i \in F^*$. For each $i \in [1,n]$, let $\tau_i$ denote the binary subform $\form{1,a_i}$ of $\psi$. By Corollary 4.6, there are rational maps $X_{\tau_i} \dashrightarrow X_\varphi$ for all $i$. In other words, $\varphi$ becomes isotropic over all the fields $F(\tau_i) \cong F(\sqrt[p]{a_i})$. By Corollary 3.2, we therefore have $a_i \in N_F(\varphi)$ for all $i$, and hence $N_F(\psi) = F^p(a_1,...,a_n) \subset N_F(\varphi)$. \end{proof} \end{proposition}

This result shows that the norm field and norm degree are birational invariants of quasilinear $p$-hypersurfaces. We will make use of this in the next section.

\section{The Izhboldin dimension and the main theorem}

In this section, we prove the main result of the paper, Theorem 5.12.\\

Let $\varphi$ be an anisotropic quasilinear $p$-form of dimension $\geq 2$ over $F$. In analogy with the theory of quadratic forms, we define the integers
\begin{equation*} \mathfrak{i}_1(X_\varphi) = \mathfrak{i}_1(\varphi) \coloneqq \mathfrak{i}_0(\varphi_{F(\varphi)}) \end{equation*}
and
\begin{equation*} \mathrm{dim}_{Izh}(X_\varphi) \coloneqq \mathrm{dim}(X_\varphi) - \mathfrak{i}_1(X_\varphi) + 1. \end{equation*}

\noindent The latter integer will be called the \emph{Izhboldin dimension} of $X_\varphi$.

\begin{example} It follows from the first part of Lemma 4.3 that $\mathfrak{i}_1(\varphi) \leq \mathrm{dim}(\varphi) - \frac{1}{p}\mathrm{dim}(\varphi)$ for any anisotropic form of dimension $\geq 2$. Generally speaking, this bound is sharp. Indeed, if $\pi$ is an anisotropic quasi-Pfister form of dimension $p^n$, then it follows from Lemma 2.6 that $\mathfrak{i}_1(\pi) = p^n - p^{n-1}$. \end{example}

\begin{example} An algebraic variety $X$ is called \emph{incompressible} if every rational morphism from $X$ to itself is dominant. In the theory of quadratic forms, an important result of A. Vishik says that any anisotropic quadric $X$ over a field of characteristic $\neq 2$ with $\mathfrak{i}_1(X) = 1$ is incompressible. This result has a key role to play in the proof of Theorem 1.3. In our setting, the corresponding statement is trivial. In fact, if $X$ is an anisotropic quasilinear $p$-hypersurface over $F$ with $\mathfrak{i}_1(X) = 1$, then $X$ only has one rational point over its function field $F(X)$. In other words, there is only one rational map from $X$ to itself. This map is, of course, the identity. \end{example}

In fact, an anisotropic quasilinear $p$-hypersurface $X$ is incompressible if and only if $\mathfrak{i}_1(X) = 1$, as the following lemma shows:

\begin{lemma} Let $\varphi$ be an anisotropic quasilinear $p$-form over $F$ of dimension $\geq 2$, and let $\psi \subset \varphi$ be a subform of codimension $ \leq \mathfrak{i}_1(\varphi) - 1$. Then the form $\psi_{F(\varphi)}$ is isotropic.
\begin{proof} The point is that the subspace of isotropic vectors for $\varphi_{F(\varphi)}$ must intersect the underlying space of $\psi_{F(\varphi)}$ non-trivially for dimension reasons. \end{proof} \end{lemma}

We are going to show (as a consequence of the main theorem) that it is impossible to find a subform of codimension larger than $\mathfrak{i}_1(\varphi) - 1$ which becomes isotropic over $F(\varphi)$. First we will need a couple of lemmas:

\begin{lemma} Let $\varphi$ be an anisotropic quasilinear $p$-form of dimension $\geq 2$ over $F$, and let $\psi \subset \varphi$ be a subform of codimension 1. Then $\varphi_{F(\varphi)} \sim \psi_{F(\varphi)}$.
\begin{proof} We can write $\varphi = \form{a} \oplus \psi$ for some $a \in F^*$. But $\psi$ represents $a$ over the field $F(\varphi)$, so the statement follows from Proposition 2.2. \end{proof} \end{lemma}

\begin{lemma} Let $\varphi$ be an anisotropic quasilinear $p$-form of dimension $\geq 2$ over $F$. Then there exists a purely transcendental field extension $K/F$ and a subform $\psi \subset \varphi_K$ of codimension $\mathfrak{i}_1(\varphi) -1$ such that $\mathfrak{i}_1(\psi) = 1$.
\begin{proof} We may assume that $\varphi = \form{1, a_1,...,a_n}$ for some $a_i \in F^*$. Let $m<n$ be such that $\mathrm{dim}(\varphi_{F(\varphi)})_{an} = m+1$. Reordering the $a_i$ if necessary, we can assume that $\form{1,a_1,...,a_m}_{F(\varphi)}$ is the anisotropic part of $\varphi_{F(\varphi)}$ by Lemma 2.4. Now, let $T_1,...,T_{n-1}$ be indeterminate variables. Then the function field $F(\varphi)$ is $F$-isomorphic to the field
\begin{equation*} F(T_1,...,T_{n-1})(\sqrt[p]{a_1T_1^p + ... + a_{n-1}T_{n-1}^p + a_n}). \end{equation*}
Let $K = F(T_{m+1},...,T_{n-1})$, and consider the subform
\begin{equation*} \psi = \form{1,a_1,...,a_m,a_{m+1}T_{m+1}^p + ... + a_{n-1}T_{n-1}^p + a_n} \end{equation*}
of $\varphi_K$. Then we have
\begin{equation*} \psi_{K(\psi)} \sim \form{1,a_1,...,a_m}_{K(\psi)} \end{equation*}
by Lemma 5.4. But the field $K(\psi)$ is $F$-isomorphic to $F(\varphi)$ by construction, so the form $\form{1,a_1,...,a_m}_{K(\psi)}$ is anisotropic. It follows that $\mathfrak{i}_1(\psi) = 1$. \end{proof} \end{lemma}

\begin{remark} We will show later (see Proposition 6.1) that $\mathfrak{i}_1(\psi) = 1$ for any subform $\psi \subset \varphi$ of codimension $\mathfrak{i}_1(\varphi) -1$ over the base field $F$. The example constructed above (modulo passing to a purely transcendental extension of $F$) will be sufficient for our more immediate concerns. \end{remark}

In the theory of quadratic forms, Theorem 1.3 is actually deduced as a consequence of the following stronger statement. It was first proved over fields of characteristic different from 2 by N. Karpenko and A. Merkurjev. It was extended to smooth quadrics in characteristic 2 in the book \cite{EKM} by R. Elman, N. Karpenko and A. Merkurjev, and to arbitrary quadrics (smooth or otherwise) by B. Totaro.

\begin{theorem}[\cite{KarpenkoMerkurjev} Theorem 3.1, \cite{Totaro1} Theorem 5.1] Let $X$ be an anisotropic quadric over a field $k$, and let $Y$ be a complete variety over $k$ which has no closed points of odd degree. Suppose that $Y$ has a closed point of odd degree over $k(X)$. Then
\begin{enumerate} \item[$\mathrm{(1)}$] $\mathrm{dim}_{Izh}(X) \leq \mathrm{dim}(Y)$.
\item[$\mathrm{(2)}$] If $\mathrm{dim}_{Izh}(X) = \mathrm{dim}(Y)$, then there exists a rational morphism $Y \dashrightarrow X$. \end{enumerate} \end{theorem}

\begin{remark} Let $X$ be an anisotropic quadric satisfying $\mathfrak{i}_1(X)=1$. In the terminology of \cite{Haution} \S 10, Theorem 5.7 says that $X$ is \emph{strongly 2-incompressible}. For an arbitrary variety $X$, the degree to which $X$ fails to satisfy the weaker notion of \emph{incompressibility} (see Example 5.2) can be measured by the minimum dimension of the image of a rational morphism from $X$ to itself. In the case where the variety $X$ is regular, this integer is commonly referred to as the \emph{canonical dimension} of $X$. The first part of Theorem 5.7 implies that the canonical dimension of a smooth anisotropic quadric is equal to the Izhboldin dimension $\mathrm{dim}_{Izh}(X)$. This includes the previously mentioned result which says that an anisotropic quadric $X$ is incompressible if and only if $\mathfrak{i}_1(X) = 1$. \end{remark}

We will now prove an analogue of Theorem 5.7 for quasilinear $p$-hypersurfaces. The approach given here is rather different from the one for quadrics given in \cite{KarpenkoMerkurjev} and \cite{Totaro1}. In fact, the latter approach does not work for quasilinear $p$-hypersurfaces whenever $p>2$. Indeed, the argument given in \cite{KarpenkoMerkurjev} and \cite{Totaro1} makes essential use of the fact that the Chow group of zero cycles on a projective quadric injects to the integers via the degree map. This is also true for isotropic quasilinear $p$-hypersurfaces, but is in general false in the anisotropic case (which is precisely the case needed for the proof). For example, if $p>2$, and $X$ is a quasilinear $p$-hypersurface of dimension 1 (i.e. a curve), then one can show that $\mathrm{CH}_0(X)$ contains nonzero $p$-torsion if and only if $X$ is regular. The proof which we present here makes no use of intersection theory, but rather exploits properties of the norm field and norm degree invariants introduced earlier. The important observation is the following lemma:

\begin{lemma} Let $\varphi$ and $\psi$ be anisotropic quasilinear $p$-forms of dimension $\geq 2$ over $F$, and let $L$ be a field such that $F \subset L \subset F(\psi)$. If the form $\varphi_{F(\psi)}$ is isotropic, then $\varphi_L$ is isotropic if and only if $\mathrm{lndeg}_L(\varphi) < \mathrm{lndeg}_F(\varphi)$.
\begin{proof} One direction was already proved in Proposition 3.1. The interesting part is the converse. So suppose that $\mathrm{lndeg}_L(\varphi) < \mathrm{lndeg}_F(\varphi)$, and suppose for the sake of contradiction that $\varphi_L$ is anisotropic. Then since $\varphi_{F(\psi)}$ is isotropic, we have $\mathrm{lndeg}_{F(\psi)}(\varphi) < \mathrm{lndeg}_L(\varphi)$ by Proposition 3.1. But this implies that $\mathrm{lndeg}_{F(\psi)}(\varphi) \leq \mathrm{lndeg}_F(\varphi) - 2$, which is impossible by the second part of Lemma 4.3. \end{proof} \end{lemma}

\begin{proposition} Let $X$ and $Y$ be anisotropic quasilinear $p$-hypersurfaces over $F$. Suppose that there exists a rational morphism $X \dashrightarrow Y$, and let $Z$ denote (the closure of) its image in $Y$. Then there exists a rational morphism $Z \dashrightarrow X$.
\begin{proof} Let $\psi$ and $\varphi$ be anisotropic forms over $F$ defining $X$ and $Y$ respectively. We need to show that $\psi$ becomes isotropic over the field $F(Z)$. But we have a natural embedding $F(Z) \subset F(\psi)$ over $F$, and since the form $\psi_{F(\psi)}$ is isotropic, Lemma 5.9 shows that it is sufficient to check that $\mathrm{lndeg}_{F(Z)}(\psi) < \mathrm{lndeg}_F(\psi)$. Now, since there exists a rational map $X \dashrightarrow Y$, we have an inclusion $N_E(\psi) \subset N_E(\varphi)$ for all extensions $E/F$ by Proposition 4.10. Suppose that $\mathrm{lndeg}_{F(Z)}(\psi) = \mathrm{lndeg}_F(\psi)$. Since $Z$ is a subvariety of $Y$, the form $\varphi_{F(Z)}$ is isotropic, and hence
\begin{equation} \label{eq1} \mathrm{lndeg}_{F(Z)}(\varphi) < \mathrm{lndeg}_F(\varphi) \end{equation}
by Proposition 3.1. By part (2) of Lemma 4.3, we also have
\begin{equation} \label{eq2} \mathrm{lndeg}_{F(\psi)}(\psi) = \mathrm{lndeg}_F(\psi) - 1 = \mathrm{lndeg}_{F(Z)}(\psi) - 1. \end{equation}
Now, since $N_{F(Z)}(\psi) \subset N_{F(Z)}(\varphi)$, it follows from \eqref{eq2} that $\mathrm{lndeg}_{F(\psi)}(\varphi) \leq \mathrm{lndeg}_{F(Z)}(\varphi) - 1$. But then combining this with \eqref{eq1}, we get $\mathrm{lndeg}_{F(\psi)}(\varphi) \leq \mathrm{lndeg}_F(\varphi) - 2$, which is impossible by the second part of Lemma 4.3. Hence $\mathrm{lndeg}_{F(Z)}(\psi) < \mathrm{lndeg}_F(\psi)$, and the proof is complete. \end{proof} \end{proposition}

\begin{corollary} Let $f \colon X \dashrightarrow Y$ be a rational morphism of anisotropic quasilinear $p$-hypersurfaces over $F$, and let $Z$ denote (the closure of) its image in $Y$. If $\mathfrak{i}_1(X) = 1$, then $X$ and $Z$ are birational via $f$.
\begin{proof} There exists a rational morphism $Z \dashrightarrow X$ by Proposition 5.10, and since $\mathfrak{i}_1(X) = 1$, the composition $X \dashrightarrow Z \dashrightarrow X$ is the identity (see Example 5.2). \end{proof} \end{corollary}

A point $y$ on a variety $Y$ over a field $k$ will be called \emph{separable} if the residue field extension $k(y)/k$ is separably generated. Here is our version of Theorem 5.7:

\begin{theorem} Let $X$ be an anisotropic quasilinear $p$-hypersurface over $F$, and let $Y$ be a complete variety over $F$ which has no separable points. Suppose that there exists a rational morphism $X \dashrightarrow Y$. Then
\begin{enumerate} \item[$\mathrm{(1)}$] $\mathrm{dim}_{Izh}(X) \leq \mathrm{dim}(Y)$.
\item[$\mathrm{(2)}$] If $\mathrm{dim}_{Izh}(X) = \mathrm{dim}(Y)$, then there exists a rational morphism $Y \dashrightarrow X$. \end{enumerate}
\begin{proof} It follows immediately from Corollary 3.5 that passing to a purely transcendental extension of the base field $F$ changes nothing in the statement. By Lemma 5.5, we may therefore assume that $X$ has a plane section $X'$ of codimension $\mathfrak{i}_1(X) - 1$ with $\mathfrak{i}_1(X') = 1$. Moreover, there exists a rational morphism $X' \dashrightarrow Y$ by Corollary 4.6. Since $\mathrm{dim}_{Izh}(X') = \mathrm{dim}_{Izh}(X)$, we may replace $X$ by $X'$ and assume that $\mathfrak{i}_1(X) = 1$ (or equivalently, $\mathrm{dim}_{Izh}(X) = \mathrm{dim}(X))$. Let $Z$ denote the image of the map $X \dashrightarrow Y$. Since $Y$ has no separable points, the field extension $F(Z)/F$ is not separably generated. It follows from Proposition 3.4 that there exists a rational morphism $Z \dashrightarrow Y'$ for some anisotropic quasilinear $p$-hypersurface $Y'$ over $F$. Now, since the given map $X \dashrightarrow Z$ is dominant, we can consider the composition $f \colon X \dashrightarrow Z \dashrightarrow Y'$. By Corollary 5.11, $X$ is birational to its image in $Y'$ via $f$. It follows that $X$ and $Z$ are birational, and part (1) of the statement follows immediately. Moreover, if we have the equality $\mathrm{dim}(X) = \mathrm{dim}(Y)$, then $X$ and $Y$ are actually birational, whence part (2). \end{proof} \end{theorem}

\begin{remark} Note that the condition that the variety $Y$ has no separable points is equivalent to the condition that $Y$ has no separable \emph{closed} points. Indeed, this follows from the well-known fact that any generically smooth variety has a dense subset of separable closed points. However, the statement of Theorem 5.12 is still weaker than the analogue of Theorem 5.7. More precisely, our result is weaker than the claim that a quasilinear $p$-hypersurface with $\mathfrak{i}_1(X)$ is \emph{strongly p-incompressible} (see Remark 5.8 and \cite{Haution} \S 10). Nonetheless, Theorem 5.12 is sufficient for several interesting applications. \end{remark}

We can now determine which subforms of an anisotropic quasilinear $p$-form $\varphi$ become isotropic over the field $F(\varphi)$. The analogue of the first part of the following statement for quadratic forms over fields of characteristic different from 2 is due to A. Vishik (see \cite{Vishik}, Corollary 3):

\begin{corollary} Let $\varphi$ be an anisotropic quasilinear $p$-form over $F$ of dimension $\geq 2$, and let $\psi$ be a subform of $\varphi$. Then $\psi_{F(\varphi)}$ is isotropic if and only if $\mathrm{codim}_{\varphi}(\psi) \leq \mathfrak{i}_1(\varphi) - 1$. In the case where $\mathrm{codim}_{\varphi}(\psi) = \mathfrak{i}_1(\varphi) - 1$, there is only one rational morphism $X_\varphi \dashrightarrow X_\psi$, and moreover, it is dominant.
\begin{proof} The fact that every subform $\psi$ of codimension $\leq \mathfrak{i}_1(\varphi) - 1$ becomes isotropic over $F(\varphi)$ was proved in Lemma 5.3. The converse follows from the first part of Theorem 5.12. If $\mathrm{codim}_{\varphi}(\psi) = \mathfrak{i}_1(\varphi) - 1$, there is only one rational morphism $X_\varphi \dashrightarrow X_\psi$, because otherwise the subspace of isotropic vectors for $\psi_{F(\varphi)}$ would have dimension $\geq 2$ and we could find a codimension $\mathfrak{i}_1(\varphi)$ subform of $\varphi$ which becomes isotropic over $F(\varphi)$. The map is dominant by part (1) of Theorem 5.12. \end{proof} \end{corollary}

We would like to prove an analogue of Theorem 1.3 for quasilinear $p$-hypersurfaces. To this end, the only obstruction is Question 4.4. We will illustrate this with a proof in the case when $p=2$ or $p=3$ (where we have a positive answer to Question 4.4 by Proposition 4.8). The case where $p=2$ was previously proved by B. Totaro in \cite{Totaro1}.

\begin{theorem} Assume that $p=2$ or $p=3$, and let $X$ and $Y$ be anisotropic quasilinear $p$-hypersurfaces over $F$. Suppose that there exists a rational morphism $X \dashrightarrow Y$. Then
\begin{enumerate} \item[$\mathrm{(1)}$] $\mathrm{dim}_{Izh}(X) \leq \mathrm{dim}_{Izh}(Y)$.
\item[$\mathrm{(2)}$] $\mathrm{dim}_{Izh}(X) = \mathrm{dim}_{Izh}(Y)$ if and only if there is a rational morphism $Y \dashrightarrow X$. \end{enumerate}
\begin{proof} Let $Y' \subset Y$ be a plane section of codimension $\mathfrak{i}_1(Y) - 1$. By Corollary 5.14, we have a dominant rational morphism $Y \dashrightarrow Y'$. Proposition 4.8 now implies that there exists a rational morphism $X \dashrightarrow Y'$. By the first part of Theorem 5.12, we get $\mathrm{dim}_{Izh}(X) \leq \mathrm{dim}(Y') = \mathrm{dim}_{Izh}(Y)$, which proves (1). If there also exists a rational morphism $Y \dashrightarrow X$, then the same argument shows that $\mathrm{dim}_{Izh}(Y) \leq \mathrm{dim}_{Izh}(X)$, and hence $\mathrm{dim}_{Izh}(X) = \mathrm{dim}_{Izh}(Y)$. On the other hand, if we are given the equality $\mathrm{dim}_{Izh}(X) = \mathrm{dim}_{Izh}(Y)$, then there is a rational morphism $Y' \dashrightarrow X$ by the second part of Theorem 5.12. Composing with the dominant rational morphism $Y \dashrightarrow Y'$, we get a rational morphism $Y \dashrightarrow X$, and this proves (2). \end{proof} \end{theorem}

\begin{remark} Note that we do not need a positive answer to Question 4.4 in its entirety to prove the analogue of Theorem 1.3. Indeed, we only need the case where $Z$ is a plane section of codimension $\mathfrak{i}_1(Y) - 1$ in $Y$. For this special case, our question is easily seen to be equivalent to:

\begin{question} Let $Y$ be an anisotropic quasilinear $p$-hypersurface over $F$. Is is true that $\mathfrak{i}_1(Y)$ is minimal among all defect indices attained by $Y$ over extensions of the base field where $Y$ becomes isotropic? \end{question} 

\noindent It is a well-established fact that $\mathfrak{i}_1(X)$ satisfies the analogous ``generic property'' when $X$ is a smooth anisotropic quadric. In our case, we have a positive answer to Question 5.17 when $p=2$ or $p=3$ by Proposition 4.8. \end{remark}

\section{Further applications to rational morphisms between quasilinear hypersurfaces}

In this section we use Theorem 5.12 to prove some more specific results concerning rational morphisms between quasilinear $p$-hypersurfaces. In particular, we prove analogues of Hoffmann's Theorem 1.1 (Theorem 6.10) and Izhboldin's Theorem 1.2 (Theorem 6.12).\\

Let $\varphi$ be an anisotropic quasilinear $p$-form over $F$ of dimension $\geq 2$, let $j \in [1,\mathfrak{i}_1(\varphi)]$, and let $\psi \subset \varphi$ be a subform of codimension $\mathfrak{i}_1(\varphi) - j$. Then it is easy to see that $\mathfrak{i}_1(\psi) \geq j$. Indeed, it suffices to show that every codimension $j-1$ subform of $\psi$ becomes isotropic over the field $F(\psi)$ (see Proposition 2.2). But any such subform is isotropic over $F(\varphi)$ by Lemma 5.3, and hence isotropic over $F(\psi)$ by Corollary 4.6. Theorem 5.12 now allows us to prove that equality holds:

\begin{proposition} Let $\varphi$ be an anisotropic quasilinear $p$-form over $F$ of dimension $\geq 2$, let $j \in [1, \mathfrak{i}_1(\varphi)]$, and let $\psi \subset \varphi$ be a subform of codimension $\mathfrak{i}_1(\varphi) - j$. Then $\mathfrak{i}_1(\psi) = j$.
\begin{proof} Suppose that $\mathfrak{i}_1(\psi) > j$. Then there is a subform $\sigma \subset \psi$ of codimension $j$ which becomes isotropic over $F(\psi)$. Let $\tau$ be a subform of $\psi$ of codimension $j-1$ which contains $\sigma$ as a codimension 1 subform. Then there exists a rational morphism $X_\tau \dashrightarrow X_\sigma$ by Corollary 4.6. Since $\tau$ has codimension $\mathfrak{i}_1(\varphi) - 1$ in $\varphi$, there is a dominant rational morphism $X_\varphi \dashrightarrow X_\tau$ by Corollary 5.14. But then taking the composition of these maps gives a rational morphism $X_\varphi \dashrightarrow X_\sigma$, and this is impossible by Corollary 5.14 \end{proof} \end{proposition}

\begin{remark} The analogue of Proposition 6.1 for quadratic forms in characteristic different from 2 was proved by A. Vishik in \cite{Vishik} prior to Karpenko and Merkurjev's Theorem 1.3. \end{remark}

We will now prove analogues of Theorems 1.1 and 1.2 for quasilinear $p$-hypersurfaces. First we need to introduce the class of \emph{quasi-Pfister neighbours}.

\begin{definition} Let $\varphi$ be a quasilinear $p$-form over $F$ of dimension $\geq 2$, and let $n$ be the unique positive integer satisfying $p^{n-1} < \mathrm{dim}(\varphi) \leq p^n$. We say that $\varphi$ is a \emph{quasi-Pfister neighbour} if $\varphi$ is proportional to a subform of a quasi-Pfister form of dimension $p^n$. In this case, the variety $X_\varphi$ will also be called a \emph{quasi-Pfister neighbour}. \end{definition}

\begin{remark} Quasi-Pfister neighbours are analogous to \emph{Pfister neighbours} in the theory of quadratic forms. For example, if $\varphi$ is a neighbour of a quasi-Pfister form $\pi$ over $F$, then for every field extension $L/F$, $\varphi_L$ is isotropic if and only if $\pi_L$ is isotropic. Indeed, if $\mathrm{dim}(\pi) = p^n$ and $\pi_L$ is isotropic, then Lemma 2.6 implies that the subspace of isotropic vectors for $\pi_L$ has dimension at least $p^n - p^{n-1}$, and therefore must intersect the underlying space of $\varphi_L$ non-trivially.  \end{remark}

Recall from Example 5.1 that if $\pi$ is an anisotropic quasi-Pfister form of dimension $p^n$, then $\mathfrak{i}_1(\pi) = p^n - p^{n-1}$. Applying Proposition 6.1 to the case of an anisotropic quasi-Pfister neighbour, we therefore get:

\begin{corollary} Let $\varphi$ be an anisotropic quasi-Pfister neighbour over $F$, and let $n$ be the unique positive integer such that $p^{n-1} < \mathrm{dim}(\varphi) \leq p^n$. Then $\mathfrak{i}_1(\varphi) = \mathrm{dim}(\varphi) - p^{n-1}$. \end{corollary}

Now, the following observation is key here:

\begin{proposition} Let $\varphi$ be an anisotropic quasilinear $p$-form of dimension $\geq 2$ over $F$. Then there exists a field extension $\widetilde{F}/F$ such that $\varphi_{\widetilde{F}}$ is an anisotropic quasi-Pfister neighbour.
\begin{proof} Let $\pi = \widehat{\nu}_F(\varphi)$ be the norm form of $\varphi$. Recall that $\varphi$ is proportional to a subform of $\pi$. Let $n$ be the unique positive integer such that $p^{n-1} < \mathrm{dim}(\varphi) \leq p^n$. If $\mathrm{dim}(\pi) = p^n$, then $\varphi$ is already a quasi-Pfister neighbour (of $\pi$). We can therefore assume that $\mathrm{dim}(\pi) \geq p^{n+1}$. In particular, we have $\mathrm{dim}_{Izh}(X_\pi) \geq p^n - 1$ by Example 5.1. Since $\mathrm{dim}(X_\varphi) \leq p^n - 2$, part (1) of Theorem 5.12 implies that there are no rational morphisms $X_\pi \dashrightarrow X_\varphi$. In other words, $\varphi$ remains anisotropic over the field $F(\pi)$. Now, the anisotropic part of $\pi_{F(\pi)}$ is nothing else but the norm form $\pi' = \widehat{\nu}_{F(\pi)}(\varphi)$ of $\varphi$ over $F(\pi)$. We have $\mathrm{dim}(\pi') < \mathrm{dim}(\pi)$, and since $\varphi_{F(\pi)}$ is anisotropic, $\varphi_{F(\pi)}$ is proportional to a subform of $\pi'$. Repeating this procedure as many times as is necessary, we eventually produce an extension $\widetilde{F}/F$ over which $\varphi$ becomes an anisotropic quasi-Pfister neighbour (of the norm form $\widehat{\nu}_{\widetilde{F}}(\varphi)$). \end{proof} \end{proposition}

\begin{remark} We should point out that the analogous statement for nondegenerate quadratic forms is certainly not true in general (see \cite{HoffmannIzhboldin} for a detailed discussion of this problem for quadratic forms over fields of characteristic different from 2). \end{remark}

\begin{corollary} Let $\varphi$ be an anisotropic quasilinear $p$-form of dimension $\geq 2$ over $F$, and let $n$ be the unique positive integer such that $p^{n-1} < \mathrm{dim}(\varphi) \leq p^n$. Then $\mathfrak{i}_1(\varphi) \leq \mathrm{dim}(\varphi) - p^{n-1}$.
\begin{proof} By Proposition 6.6, there is an extension $\widetilde{F}/F$ such that $\varphi_{\widetilde{F}}$ is an anisotropic quasi-Pfister neighbour. By Corollary 6.5, we therefore have $\mathfrak{i}_1(\varphi_{\widetilde{F}}) = \mathrm{dim}(\varphi) - p^{n-1}$. Since we have a natural embedding $F(\varphi) \subset \widetilde{F}(\varphi_{\widetilde{F}})$, we get
\begin{equation*} \mathfrak{i}_1(\varphi) = \mathfrak{i}_0(\varphi_{F(\varphi)}) \leq 
\mathfrak{i}_0(\varphi_{\widetilde{F}(\varphi_{\widetilde{F}})}) = \mathfrak{i}_1(\varphi_{\widetilde{F}}) = \mathrm{dim}(\varphi) - p^{n-1}, \end{equation*}
which is what we wanted. \end{proof} \end{corollary}

\begin{example} Let $\varphi$ be an anisotropic form of dimension $p^n + 1$ for some $n \geq 0$. Then $\mathfrak{i}_1(\varphi) = 1$. \end{example}

Now we can prove an analogue of Theorem 1.1 for quasilinear $p$-hypersurfaces:

\begin{theorem} Let $X$ and $Y$ be anisotropic quasilinear $p$-hypersurfaces over $F$. If there exists $n \geq 1$ such that $\mathrm{dim}(Y) \leq p^n - 2 < \mathrm{dim}(X)$, then there are no rational morphisms $X \dashrightarrow Y$. 
\begin{proof} By part (1) of Theorem 5.12, it is sufficient to show that $\mathrm{dim}_{Izh}(X) > p^n - 2$, and this follows from Corollary 6.8. \end{proof} \end{theorem}

\begin{remark} The analogue of Corollary 6.8 for quadratic forms over fields of characteristic different from 2 was originally proved by D. Hoffmann as a corollary of Theorem 1.1 (see \cite{Hoffmann1}). Here the roles are reversed. The difference is that we were able to use Theorem 5.12 to prove Corollary 6.6, but, as we have remarked above, the analogue of Corollary 6.6 for nondegenerate quadratic forms is false in general. Corollary 6.8 and Theorem 6.10 were proved in the case $p=2$ by D. Hoffmann and A. Laghribi in \cite{HoffmannLaghribi2} using different methods. \end{remark}

We also get the following analogue of Izhboldin's Theorem 1.2. The case where $p=2$ was proved using different methods by D. Hoffmann and A. Laghribi in \cite{HoffmannLaghribi2}.

\begin{theorem} Let $X$ and $Y$ be anisotropic quasilinear $p$-hypersurfaces over $F$ with $\mathrm{dim}(Y) = p^n - 1 \leq \mathrm{dim}(X)$ for some $n \geq 0$. Suppose that there exists a rational morphism $X \dashrightarrow Y$. Then there exists a rational morphism $Y \dashrightarrow X$. If in addition we have $\mathrm{dim}(X) = p^n -1$, then $X$ and $Y$ are birational.
\begin{proof} Let $X' \subset X$ be a plane section of dimension $p^n - 1$. By Corollary 4.6, there exists a rational morphism $X' \dashrightarrow Y$. By Example 6.9, we have that $\mathfrak{i}_1(X') = 1$. It then follows from Corollary 5.11 that $X'$ is birational to $Y$, whence the result. \end{proof} \end{theorem}

Note that in the case where $\mathrm{dim}(X) = \mathrm{dim}(Y) = p^n -1$ for some $n \geq 0$, we get the stronger assertion (in comparison with Theorem 1.2) that $X$ and $Y$ are birational. For non-quasilinear quadrics, this is still an open problem (see also Conjecture 7.1). 

\begin{remark} In a recent article \cite{Haution}, O. Haution has shown that any degree $p$ hypersurface of dimension $p-1$ (over an arbitrary field) which has no closed points of degree prime to $p$ is strongly $p$-incompressible. In particular, an anisotropic quasilinear $p$-hypersurface of dimension $p-1$ is strongly $p$-incompressible. For such varieties, this statement is stronger than the results established above in the present article (see Remark 5.13). The proof uses an extension of K. Zainoulline's degree formula for the Euler characteristic (see \cite{Zainoulline}) to fields of arbitrary characteristic, and also to non-smooth varieties of sufficiently small dimension. The Euler characteristic (of the structure sheaf) does not, however, distinguish between degree $p$ hypersurfaces of dimension $\geq p-1$, so this invariant can only be used to explain the first level of the ``separation'' exhibited by Theorems 6.10 and 6.12. We remark any smooth degree $p$ hypersurface of dimension $p^n -1$ (for any $n\geq 0$) over a field of characteristic $\neq p$ which has no closed points of degree prime to $p$ is strongly $p$-incompressible by the degree formulas of A. Merkurjev (generalising those of Rost; see \cite{Merkurjev}). It is not known at present if the degree formulas used to prove this hold in arbitrary characteristic. \end{remark}

\section{Birational geometry of quasilinear hypersurfaces}

In this section we consider the extension of the results obtained by B. Totaro on the birational geometry of quadrics in \cite{Totaro1} to quasilinear hypersurfaces of higher degree.\\

An old problem of O. Zariski asks whether two stably birational varieties of the same dimension over a field are actually birational. It is well-known that this is false in general, but in the case where both varieties are smooth anisotropic quadrics, it is still an important open problem. Using the fact that a smooth isotropic quadric is a rational variety, one easily shows that two smooth anisotropic quadrics $X$ and $Y$ over a field are stably birational if and only if there exist rational morphisms $X \dashrightarrow Y$ and $Y \dashrightarrow X$. We can therefore ask the following question for arbitrary quadrics:

\begin{conjecture}[Quadratic Zariski Problem] Let $X$ and $Y$ be anisotropic quadrics of the same dimension over a field $k$. Suppose that there exist rational morphisms $X \dashrightarrow Y$ and $Y \dashrightarrow X$. Is it true that $X$ and $Y$ are birational? \end{conjecture}

In a series of papers (\cite{Totaro3}, \cite{Totaro1}, \cite{Totaro2}), B. Totaro has suggested a new approach to this problem by means of a related conjecture concerning rulings on quadrics. We will say that a variety $X$ over a field $k$ is \emph{ruled} if $X$ is birational to $Y \times \mathbb{P}_k^1$ for some variety $Y$ over $k$. Totaro has observed the following consequence of Vishik's result stating that a quadric with first Witt index equal to 1 is incompressible:

\begin{theorem}[\cite{Totaro1}, Corollary 3.2] Let $X$ be an anisotropic quadric over a field $k$. If $\mathfrak{i}_1(X) = 1$, then $X$ is not ruled. \end{theorem}

Moreover he conjectures:

\begin{conjecture}[\cite{Totaro3}, Conjecture 3.1] Let $X$ be an anisotropic quadric over a field $k$. Then $X$ is ruled if and only if $\mathfrak{i}_1(X) > 1$. \end{conjecture}

This is formulated more precisely as follows:

\begin{conjecture}[\cite{Totaro2}, Conjecture 1.1] Let $X$ be an anisotropic quadric over a field $k$. Then $X$ is birational to $X' \times \mathbb{P}_k^{\mathfrak{i}_1(X) - 1}$ for some subquadric $X' \subset X$ of codimension $\mathfrak{i}_1(X) - 1$. \end{conjecture}

Many results are known on all these problems for quadrics of small dimension. In \cite{Totaro1}, Totaro proves Conjecture 7.4 for the entire class of quasilinear quadrics, and then uses it to prove the quasilinear case of Conjecture 7.1. The same approach should work for quasilinear hypersurfaces of higher degree. The analogue of Theorem 7.2 is trivial in this setting:

\begin{proposition} Let $X$ be an anisotropic quasilinear $p$-hypersurface over $F$. If $\mathfrak{i}_1(X) = 1$, then $X$ is not ruled.
\begin{proof} Suppose that $X$ is birational to $Y \times \mathbb{P}_F^1$ for some variety $Y$ over $F$. Then we can construct a rational morphism $X \dashrightarrow X$ as the composition
\begin{equation*} X \dashrightarrow Y \times \mathbb{P}_F^1 \xrightarrow{pr_Y} Y \hookrightarrow Y \times \mathbb{P}_F^1 \dashrightarrow X, \end{equation*}
where the second map is the canonical projection and the third map is the embedding of $Y$ in $Y \times \mathbb{P}_F^1$ at a rational point of $\mathbb{P}_F^1$. Note that the composition is defined because the third map embeds $Y$ as an effective Cartier divisor in $Y \times \mathbb{P}_F^1$ (see Proposition 4.5). But the resulting map is not surjective by construction, and this contradicts the fact that the only rational morphism from $X$ to itself is the identity (see Example 5.2). \end{proof} \end{proposition}

We can now prove an analogue of Conjecture 7.4 for quasilinear hypersurfaces of higher degree. Given the results of \S 5, the proof for the case $p=2$ given by B. Totaro in \cite{Totaro1} carries over verbatim to all primes $p$. We reproduce the argument here for the reader's convenience:

\begin{theorem} Let $X$ be an anisotropic quasilinear $p$-hypersurface over $F$, and let $X' \subset X$ be a plane section of codimension $\mathfrak{i}_1(X) - 1$. Then $X$ is birational to $X' \times \mathbb{P}_F^{\mathfrak{i}_1(X) - 1}$.
\begin{proof} For simplicity of notation, let us put $r = \mathfrak{i}_1(X)$. Let $\varphi$ be an anisotropic form over $F$ which defines the variety $X$, and let $\psi \subset \varphi$ be a subform of codimension $r-1$ which defines its subvariety $X'$. By Corollary 5.14, there exists a dominant rational morphism $\pi \colon X \dashrightarrow X'$ via which we may view $F(\psi)$ as a subfield of $F(\varphi)$. In particular, the defect index of $\varphi$ over $F(\psi)$ is no more than $r$. On the other hand, it is at least this large by a simple application of Corollary 4.6. Hence $\varphi$ has the same defect index over $F(\psi)$ as it does over $F(\varphi)$. Let $v_0,...,v_{r-1}$ be a basis of the subspace of isotropic vectors for the form $\varphi_{F(\psi)}$. The $v_i$ may be regarded as rational morphisms from $X'$ to the affine hypersurface $\lbrace \varphi = 0 \rbrace$. Define a rational morphism $f \colon X' \times \mathbb{P}_F^{r-1} \dashrightarrow X$ over $F$ by the assignment
\begin{equation*} (x',[\lambda_0:...:\lambda_{r-1}]) \mapsto [\lambda_0 v_0(x') + ... + \lambda_{r-1} v_{r-1}(x')]. \end{equation*}
Now, the identity map $X \rightarrow X$ corresponds to some isotropic line in the space $V_\varphi \otimes _F F(\varphi)$. Since $\varphi$ has the same index over $F(\psi)$ as it does over $F(\varphi)$, there are rational functions $f_i \in F(\varphi)$ such that $[f_0 \cdot (v_0 \circ \pi) + ... + f_{r-1} \cdot (v_{r-1} \circ \pi)]$ is the identity map from $X$ to itself. Define a rational morphism $g \colon X \dashrightarrow X' \times \mathbb{P}_F^{r-1}$ by the assignment
\begin{equation*} x \mapsto (\pi(x),[f_0(x),...,f_{r-1}(x)]). \end{equation*}
By construction, the composition $f \circ g$ is the identity on $X$. Therefore $g$ is a birational isomorphism, and the statement is proved. \end{proof} \end{theorem}

Finally, we make some remarks concerning an analogue of the Quadratic Zariski Problem for quasilinear hypersurfaces of higher degree. Unfortunately, the obstruction here is again Question 4.4. As with Theorem 5.15, we illustrate this with a proof for quasilinear quadrics and cubics (for which we know that Question 4.4 has a positive answer).

\begin{theorem} Assume that $p=2$ or $p=3$, and let $X$ and $Y$ be anisotropic quasilinear $p$-hypersurfaces of the same dimension over $F$. Suppose that there exist rational morphisms $X \dashrightarrow Y$ and $Y \dashrightarrow X$. Then $X$ and $Y$ are birational.
\begin{proof} Let $X' \subset X$ and $Y' \subset Y$ be plane sections of codimensions $\mathfrak{i}_1(X) - 1$ and $\mathfrak{i}_1(Y) - 1$ respectively. By Theorem 7.6, $X$ is birational to $X' \times \mathbb{P}_F^{\mathfrak{i}_1(X) - 1}$ and $Y$ is birational to $Y' \times \mathbb{P}_F^{\mathfrak{i}_1(Y) - 1}$. Moreover, we have $\mathfrak{i}_1(X) = \mathfrak{i}_1(Y)$ by Theorem 5.15 (here we are using the statement of Question 4.4). In order to prove the statement, it therefore suffices to show that $X'$ and $Y'$ are birational. Now, we have rational morphisms $X' \dashrightarrow Y$ and $Y \dashrightarrow Y'$ by Corollary 4.6 and Lemma 5.3 respectively. By Proposition 4.8, we therefore have a rational morphism $X' \dashrightarrow Y'$ (again, we are using the statement of Question 4.4). Since $\mathrm{i}_1(X) = 1$ and $\mathrm{dim}(X') = \mathrm{dim}(Y')$, $X'$ and $Y'$ are birational by Corollary 5.11. \end{proof} \end{theorem}

\begin{remark} As with Theorem 5.15, we do not need a positive answer to Question 4.4 in its entirety, but only a positive answer to Question 5.17. \end{remark}

Still, using some of the results we have obtained, we can give partial results towards a positive solution to the Zariski problem for all primes $p$. By Corollary 5.11, the conjecture is true whenever $\mathfrak{i}_1(X) = 1$. In particular, it is true whenever $\mathrm{dim}(X) = \mathrm{dim}(Y) = p^n - 1$ for some $n \geq 0$ by Example 6.9. We can improve this to include dimensions which are ``sufficiently close'' to the form $p^n - 1$. First we need the case of quasi-Pfister neighbours:

\begin{proposition} Let $X$ and $Y$ be quasi-Pfister neighbours of the same dimension over $F$. Suppose that there are rational morphisms $X \dashrightarrow Y$ and $Y \dashrightarrow X$. Then $X$ and $Y$ are birational.
\begin{proof} Let $X' \subset X$ and $Y' \subset Y$ be plane sections of codimensions $\mathfrak{i}_1(X) - 1$ and $\mathfrak{i}_1(Y) - 1$ respectively. From the proof of Theorem 7.7, the only thing left to check is that we have rational morphisms $X' \dashrightarrow Y'$ and $Y' \dashrightarrow X'$. Now, let $\varphi$ and $\psi$ be quasi-Pfister forms over $F$ such that $X$ is a neighbour of $X_{\varphi}$ and $Y$ is a neighbour of $X_{\psi}$. Then there exist rational morphisms $X' \dashrightarrow X_{\psi}$ and $Y' \dashrightarrow X_{\varphi}$ by Corollary 4.6. But by Corollary 6.5, $X'$ and $Y'$ are also neighbours of $X_{\varphi}$ and $X_{\psi}$ respectively. In particular, for any variety $Z$ over $F$, there exists a rational morphism $Z \dashrightarrow X'$ (resp. $Z \dashrightarrow Y'$) if and only if there exists a rational morphism $Z \dashrightarrow X_{\varphi}$ (resp. $Z \dashrightarrow X_{\psi}$; see Remark 6.4). Therefore we have rational morphisms $X' \dashrightarrow Y'$ and $Y' \dashrightarrow X'$, and the proof is complete. \end{proof} \end{proposition}

\begin{remark} For example, any two neighbours of the same quasi-Pfister hypersurface which have the same dimension are birational. Together with Proposition 6.6, Proposition 7.9 shows that the analogue of the quadratic Zariski problem is true up to making an extension of the base field which preserves the anisotropy of the quasilinear hypersurfaces involved. \end{remark}

We conclude with the following result, which settles the Zariski problem in a large number of cases:

\begin{proposition} Let $X$ and $Y$ be anisotropic quasilinear $p$-hypersurfaces of the same dimension $d$ over $F$. Suppose that there exist rational morphisms $X \dashrightarrow Y$ and $Y \dashrightarrow X$, and let $n$ be the unique non-negative integer satisfying $p^n < d+2 \leq p^{n+1}$. If $d \leq p^n + n$, then $X$ and $Y$ are birational.
\begin{proof} Let $\varphi$ and $\psi$ be anisotropic forms defining $X$ and $Y$ respectively. By Proposition 4.10, we have $\mathrm{lndeg}_F(\varphi) = \mathrm{lndeg}_F(\psi)$. Let us denote this integer by $m$. If $m = {n+1}$, then $X$ and $Y$ are quasi-Pfister neighbours, and we are done by Proposition 7.9. We can therefore assume that $m \geq n+2$. Now, let $X' \subset X$ and $Y' \subset Y$ be plane sections of codimensions $\mathfrak{i}_1(X) - 1$ and $\mathfrak{i}_1(Y) - 1$ respectively. As before, we only need to show that there exist rational morphisms $X' \dashrightarrow Y'$ and $Y' \dashrightarrow X'$. By Corollary 4.6 and Lemma 5.3, we have rational morphisms $X' \dashrightarrow Y$ and $Y \dashrightarrow Y'$ (resp. rational morphisms $Y' \dashrightarrow X$ and $X \dashrightarrow X'$), and by Proposition 4.5 it will be sufficient to prove that $Y$ (resp. $X$) is regular at the image of the generic point of $X'$ (resp. $Y'$). We may therefore assume that $X$ and $Y$ are not regular. Note that we have $\mathfrak{i}_1(X') = \mathfrak{i}_1(Y') = 1$ by Proposition 6.1. It therefore follows from Corollary 5.11 that we may also assume that the given rational morphisms $X' \dashrightarrow Y$ and $Y' \dashrightarrow X$ are closed embeddings of subvarieties. Now, the non-regular locus of $X$ (resp. $Y$) has codimension $m \geq n+2$ by Theorem 4.2. On the other hand, the Separation Theorem 6.10 implies that
\begin{equation*} \mathrm{codim}_Y(X') = d - \mathrm{dim}(X') \leq d - (p^n - 1), \end{equation*}
and similarly $\mathrm{codim}_X(Y') \leq d - (p^n - 1)$. Hence if $d \leq p^n + n$, then we have
\begin{equation*} \mathrm{codim}_Y(X'), \mathrm{codim}_X(Y') \leq n+1 < m, \end{equation*}
and so $Y$ (resp. $X$) is regular at the generic point of $X'$ (resp. $Y'$). \end{proof} \end{proposition}

{\bf Acknowledgements.} This work is part of my PhD thesis carried out at the University of Nottingham. I would like to thank Detlev Hoffmann for suggesting several of the problems discussed in this text and for his guidance during the period in which this work was completed. I would also like to thank Olivier Haution and Alexander Vishik for numerous helpful discussions. I am particularly grateful to Alexander Vishik for comments which helped to improve \S 5 of the paper.\\

\bibliographystyle{alphaurl}
\bibliography{RMBQHA}

\begin{thebibliography}{EKM08}

\bibitem[BK86]{BlochKato}
S.~Bloch and K.~Kato.
\newblock {$p$-adic \'{e}tale cohomology}.
\newblock {\em Inst. Hautes \'{E}tudes Sci. Publ. Math.}, (63):107--152, 1986.

\bibitem[EKM08]{EKM}
R.~Elman, N.~Karpenko, and A.~Merkurjev.
\newblock {\em {The algebraic and geometric theory of quadratic forms}}.
\newblock Number~56 in American Mathematical Society Colloquium Publications.
  American Mathematical Society, 2008.

\bibitem[Hau11]{Haution}
O.~Haution.
\newblock {Integrality of the Chern character in small codimension}.
\newblock Preprint, 2011.
\newblock \href {http://arxiv.org/abs/1103.4084v3} {\path{arXiv:1103.4084v3}}.

\bibitem[HI00]{HoffmannIzhboldin}
D.W. Hoffmann and O.T. Izhboldin.
\newblock {Embeddability of quadratic forms in Pfister forms}.
\newblock {\em Indag. Math. (N.S.)}, 11(2):219--237, 2000.

\bibitem[HL04]{HoffmannLaghribi1}
D.W. Hoffmann and A.~Laghribi.
\newblock {Quadratic forms and Pfister neighbors in characteristic 2}.
\newblock {\em Trans. Amer. Math. Soc.}, 356(10):4019--4053, 2004.

\bibitem[HL06]{HoffmannLaghribi2}
D.W. Hoffmann and A.~Laghribi.
\newblock {Isotropy of quadratic forms over the function field of a quadric in
  characteristic 2}.
\newblock {\em J. Algebra}, 295(2):362--386, 2006.

\bibitem[Hof95]{Hoffmann1}
D.W. Hoffmann.
\newblock {Isotropy of quadratic forms over the function field of a quadric}.
\newblock {\em Math. Z.}, 220(3):461--476, 1995.

\bibitem[Hof04]{Hoffmann2}
D.W. Hoffmann.
\newblock {Diagonal forms of degree $p$ in characteristic $p$}.
\newblock In {\em Algebraic and arithmetic theory of quadratic forms}, Contemp.
  Math., {\bf 344}, pages 135--183. Amer. Math. Soc., 2004.

\bibitem[Izh00]{Izhboldin}
O.T. Izhboldin.
\newblock {Motivic equivalence of quadratic forms {II}}.
\newblock {\em Manuscripta Math.}, 102(1):41--52, 2000.

\bibitem[Kar03]{Karpenko}
N.A. Karpenko.
\newblock {On the first Witt index of quadratic forms}.
\newblock {\em Invent. Math.}, 153(2):455--462, 2003.

\bibitem[Kat82]{Kato}
K.~Kato.
\newblock {Galois cohomology of complete discrete valuation fields}.
\newblock In {\em Algebraic $K$-theory, Part II (Oberwolfach, 1980)}, Lecture
  Notes in Math., pages 215--238. Springer, 1982.

\bibitem[KM03]{KarpenkoMerkurjev}
N.~Karpenko and A.~Merkurjev.
\newblock {Essential dimension of quadrics}.
\newblock {\em Invent. Math.}, 153(2):361--372, 2003.

\bibitem[Lan02]{Lang}
S.~Lang.
\newblock {\em {Algebra}}.
\newblock Number 211 in Graduate Texts in Mathematics. Springer Verlag, 2002.
\newblock Revised Third Edition.

\bibitem[LM07]{LevineMorel}
M.~Levine and F.~Morel.
\newblock {\em {Algebraic cobordism}}.
\newblock Springer Monographs in Mathematics. Springer, 2007.

\bibitem[Mer00]{Rost}
A.~Merkurjev.
\newblock {Degree formula}.
\newblock Preprint, 2000.

\bibitem[Mer03]{Merkurjev}
A.~Merkurjev.
\newblock {Steenrod operations and degree formulas}.
\newblock {\em J. Reine Angew. Math.}, 565:13--26, 2003.

\bibitem[Sch10]{Schroer}
S.~Schr$\ddot{\mathrm{o}}$er.
\newblock {On fibrations whose geometric fibers are nonreduced}.
\newblock {\em Nagoya Math. J.}, 200:35--57, 2010.

\bibitem[Tot07]{Totaro3}
B.~Totaro.
\newblock {The automorphism group of an affine quadric}.
\newblock {\em Math. Proc. Cambridge Philos. Soc.}, 143(1):1--8, 2007.

\bibitem[Tot08]{Totaro1}
B.~Totaro.
\newblock {Birational geometry of quadrics in characteristic 2}.
\newblock {\em J. Algebraic Geom.}, 17(3):577--597, 2008.

\bibitem[Tot09]{Totaro2}
B.~Totaro.
\newblock {Birational geometry of quadrics}.
\newblock {\em Bull. Soc. Math. France}, 137(2):253--276, 2009.

\bibitem[Vis99]{Vishik}
A.~Vishik.
\newblock {Direct summands in the motives of quadrics}.
\newblock Preprint, 1999.

\bibitem[Zai10]{Zainoulline}
K.~Zainoulline.
\newblock {Degree formula for connective $K$-theory}.
\newblock {\em Invent. Math.}, 179(3):507--522, 2010.

\end{thebibliography}

\end{document}